\documentclass[12pt]{iopart}
\usepackage{amsfonts}
\usepackage[english]{babel}
\usepackage{euscript}
\newcommand{\textfrac}[2]{{\textstyle\frac{#1}{#2}}}
\begin{document}

\title[Symmetries and coverings for the r-th mdKP equation]{%
Contact Integrable Extensions of Symmetry Pseudo-Group and Coverings of the r-th Modified Dispersionless Kadomtsev -- Petviashvili Equation
}

\author{Oleg I. Morozov}

\address{Department of Mathematics, Moscow State Technical University of Civil Aviation, Kronshtadtskiy Blvd 20, Moscow 125993, Russia
\\
oim{\symbol{64}}foxcub.org}

\begin{abstract}
We apply the technique of integrable extensions to the symmetry pseudo-group of the r-th mdKP equation. This gives another look on deriving known coverings and allows us to find new coverings for this equation.
%\keywords{Lie pseudo-groups \and Maurer--Cartan forms \and symmetries of differential equations 
%\and coverings of differential equations}
\end{abstract}


\maketitle

\section{Introduction}

Many useful methods to study nonlinear partial differential equations ({\sc pde}s) such as  
inverse scattering transformations, B\"ack\-lund transformations, recursion operators, nonlocal symmetries and non\-lo\-cal conservation laws, can be successfully described in the framework of the theory of coverings, \cite{KV84,KLV,KV89,KV99}. Consequently, a problem to find a covering for a given {\sc pde} is of a significant importance. A number of techniques have been devised to handle this problem,  
\cite{WE,Estabrook,DoddFordy,Hoenselaers,Marvan1997,Marvan2002,Sakovich,Igonin,%
Morris1976,Morris1979,Zakharov82,Tondo,Marvan1992,Harrison1995,Harrison2000}.
In \cite{Kuzmina}, examples of coverings of {\sc pde}s with three independent variables were found by means of \'Elie Cartan's method of equivalence, \cite{Cartan1,Cartan2,Cartan3,Cartan4,Gardner,Kamran,Olver95}. This approach was developed in \cite{Morozov2007,Morozov2008,Morozov2008b}, where it was shown that for a number of {\sc pde}s their coverings can be inferred from invariant linear combinations of Maurer--Cartan  ({\sc mc}) forms  of the contact symmetry pseudo-groups of the {\sc pde}s. The aim of the present paper is to make the method of \cite{Morozov2007,Morozov2008,Morozov2008b} more systematic. We apply the technique of integrable extensions of exterior differential systems, \cite[\S 6]{BryantGriffiths}, to {\sc mc} forms of the symmetry pseudo-group of the r-th modified dispersionless Kadomtsev--Petviashvili equation, \cite{Blaszak}, and derive the known coverings, \cite{FerapontovKhusnutdinova,ChangTu,KonopelchenkoAlonso,Dunajski,Morozov2008b}, and three new coverings.

\section{Preliminaries}

\subsection{Coverings of PDEs}

Let $\pi_{\infty} : J^{\infty}(\pi) \rightarrow \mathbb{R}^n$ be the infinite jet bundle of local sections of the bundle $\pi : \mathbb{R}^n \times \mathbb{R} \rightarrow \mathbb{R}$. The coordinates on $J^{\infty}(\pi)$ are $(x^i, u_I)$, where
$I = (i_1,...,i_k)$ are symmetric multi-indices,  $i_1,...,i_k \in \{1,...,n\}$, $u_\emptyset= u$, and for any local section $f$ of $\pi$ there exists a section $j_{\infty}(f) : \mathbb{R}^n \rightarrow J^{\infty}(\pi)$ such that $u_I(j_{\infty}(f)) = \partial^{\# I}(f)/\partial x^{i_1} ... \partial x^{i_k}$, $\# I =\#(i_1,...,i_k) = k$.  The {\it total derivatives} on $J^{\infty}(\pi)$ are defined in the local coordintes as
\[
D_i = \frac{\partial}{\partial x^i}
+\sum \limits_{\# I \ge 0} u_{Ii} \, \frac{\partial}{\partial u_I}.
\]
\noindent
We have $[D_i, D_j] = 0$ for $i, j \in \{1,...,n\}$. A differential equation $F(x^i, u_K)=0$  defines a submanifold 
$\EuScript{E}^{\infty} = \{ D_I(F) =0 \,\,\vert\,\, \#I\ge 0\} \subset J^{\infty}(\pi)$,
where $D_I = D_{i_1}\circ ... \circ D_{i_k}$ for $I=(i_1,...,i_k)$. We denote restrictions of $D_i$ on $\EuScript{E}^{\infty}$ as $\bar{D}_i$.

In local coordinates, a {\it covering} over $\EuScript{E}^{\infty}$ is a bundle
$\widetilde{\EuScript{E}}^{\infty} = \EuScript{E}^{\infty} \times \EuScript{V} \rightarrow \EuScript{E}^{\infty}$ with fibre coordinates $v^\alpha$, $\alpha \in \{1,..., N\}$ or $\alpha \in \mathbb{N}$, equipped with {\it ex\-ten\-ded total derivatives}
\[
\widetilde{D}_i = \bar{D}_i
+\sum \limits_{\alpha}
T^{\alpha}_i (x^j, u_I, v^{\beta})\,\frac{\partial}{\partial v^\alpha}
\]
such that $[\widetilde{D}_i, \widetilde{D}_j ]=0$ whenever $(x^i, u_I) \in \EuScript{E}^{\infty}$.

In terms of differential forms, the covering is defined by the forms, \cite{WE},
\[
\omega^{\alpha} = d v^{\alpha}- T^{\alpha}_i (x^j, u_I, v^{\beta})\,dx^i
\]
\noindent
such that
\begin{equation}
d \omega^{\alpha} \equiv 0 \,\,\,{\mathrm{mod}}\,\,\,\omega^{\beta}, \bar{\vartheta}_I 
\quad \Longleftrightarrow \quad (x^i,u_I) \in \EuScript{E}^{\infty},
\label{WE_def}
\end{equation}
where $\bar{\vartheta}_I$ are restrictions of contact forms 
$\vartheta_I = du_I-u_{I,k}\,dx^k$ on $\EuScript{E}^{\infty}$.

\subsection{Cartan's structure theory of Lie pseudo-groups}

Let $M$ be a manifold of dimension $n$. A {\it local diffeomorphism} on $M$ is a diffeomorphism 
$\Phi : \EuScript{U} \rightarrow \widetilde{\EuScript{U}}$ of two open subsets of $M$. A {\it pseudo-group} $\mathfrak{G}$ on $M$ is a collection of local dif\-feo\-mor\-phisms of $M$, which is closed under composition {\it when defined}, contains an identity and is closed under inverse. A {\it Lie pseudo-group} is a pseudo-group whose diffeomorphisms are local analytic solutions of an involutive system of partial differential equations. 

\'Elie Cartan's approach to Lie pseudo-groups is based on a possibility to characterize transformations from a pseudo-group in terms of a set of invariant differential 1-forms called {\it Maurer--Cartan forms}. In a general case {\sc mc} forms $\omega^1$, ..., $\omega^m$ of a Lie pseudo-group $\mathfrak{G}$ on $M$ are defined on a direct product $\hat{M} \times G$, where $\mu:\hat{M}\times G \rightarrow M$ is a bundle, $m = \dim \,\hat{M}$, $G$ is a finite-dimensional Lie group. The forms $\omega^i$ are semi-basic w.r.t. the natural projection 
$\hat{M} \times G \rightarrow \hat{M}$ and define a coframe on $\hat{M}$, that is, a basis of 
the cotangent bundle of $\hat{M}$. They characterize the pseudo-group $\mathfrak{G}$ in the following sense: a local diffeomorphism 
$\Phi : \EuScript{U} \rightarrow \widetilde{\EuScript{U}}$ on $M$ belongs to $\mathfrak{G}$ whenever there exists a local diffeomorphism 
$\Psi : \EuScript{V} \rightarrow \widetilde{\EuScript{V}}$ on $\hat{M}\times G$ such that 
$\mu \circ \Psi = \Phi \circ \mu$ and the forms $\omega^j$ are invariant w.r.t. $\Psi$, that is,
\begin{equation}
\Psi^{*} \left(\omega^i\vert {}_{\widetilde{\EuScript{V}}} \right) 
= \omega^i\vert {}_{\EuScript{V}}.
\label{Phi_omega}
\end{equation}
Expressions of exterior differentials of the forms $\omega^i$ in terms of themselves give Cartan's structure equations of $\mathfrak{G}$:
\begin{equation}
d \omega^i = A_{\alpha j}^i \,\pi^\alpha \wedge \omega^j + B_{jk}^i\,\omega^j \wedge \omega^k,
\qquad B_{jk}^i = - B_{kj}^i.
\label{SE} 
\end{equation}
Here and below we assume summation on repeated indices. The forms $\pi^\alpha$, 
$\alpha \in \{1,...,\dim \, G\}$, are linear combinations of {\sc mc} forms of the Lie group $G$ and the forms $\omega^i$. The coefficients $A_{\alpha j}^i$ and $B_{jk}^i$ are either constants or functions of a set of invariants $U^\kappa : M \rightarrow \mathbb{R}$, 
$\kappa \in \{1,...,l\}$, $l < \dim\, M$, of the pseudo-group $\mathfrak{G}$, so 
$\Phi^{*} \left(U^{\kappa}\vert {}_{\widetilde{\EuScript{U}}} \right) 
= U^{\kappa}\vert {}_{\EuScript{U}}$
for every $\Phi \in \mathfrak{G}$. In the latter case, differentials of $U^\kappa$ are invariant 1-forms, so they are linear combinations of the forms $\omega^j$,
\begin{equation}
d U^\kappa = C_j^\kappa\,\omega^j,
\label{dUs}
\end{equation} 
where the coefficients $C_j^\kappa$ depend on the invariants $U^1$, ..., $U^l$ only.

Eqs. (\ref{SE}) must be compatible in the following sense: we have  
\begin{equation}
d(d \omega^i) = 0 = d \left(A_{\alpha j}^i \,\pi^\alpha \wedge \omega^j + B_{jk}^i\,\omega^j \wedge \omega^k \right),
\label{compatibility_conditions_SE} 
\end{equation}
therefore there must exist expressions 
\begin{equation}
d \pi^\alpha = W_{\lambda j}^\alpha\, \chi^\lambda \wedge \omega^j 
+ X_{\beta \gamma}^\alpha\,\pi^\beta\wedge\pi^\gamma
+Y_{\beta j}^\alpha\,\pi^\beta \wedge \omega^j
+Z_{jk}^\alpha\,\omega^j \wedge \omega^k,
\label{prolonged_SE}
\end{equation}
with some additional 1-forms $\chi^\lambda$ and the coefficients $W_{\lambda j}^\alpha$ to 
$Z_{j k}^\alpha$ depending on the invariants $U^\kappa$ such that the right-hand side of (\ref{compatibility_conditions_SE}) appear to be identically equal to zero after substituting for (\ref{SE}), (\ref{dUs}), and (\ref{prolonged_SE}). Also, from (\ref{dUs}) it follows that
the right-hand side of the equation
\begin{equation}
d(d U^\kappa) = 0 = d(C_j^\kappa\,\omega^j)
\label{compatibility_conditions_dUs}
\end{equation}
must be identically equal to zero after substituting for (\ref{SE}) and (\ref{dUs}).

The forms $\pi^\alpha$ are not invariant w.r.t. the pseudo-group $\mathfrak{G}$. Respectively, the structure equations (\ref{SE}) are not changing when replacing 
$\pi^\alpha \mapsto \pi^\alpha + z^\alpha_j\,\omega^j$ for certain parametric coefficients $z^\alpha_j$. The dimension $r^{(1)}$ of the linear space of these coefficients satisfies the fol\-lo\-wing inequality
\begin{equation}
r^{(1)} \le n\,\dim \, G  -  \sum \limits_{k=1}^{n-1} (n-k)\,\sigma_k, 
\label{CT}
\end{equation}
where the {\it reduced characters} $\sigma_k$ are defined by
\[
\sigma_k = \max \limits_{u_1,...,u_k}\, \mathrm{rank}\,\, \mathbb{A}_k(u_1,...,u_k)
\]
with the  matrices $\mathbb{A}_k$ inductively defined by
\[
\mathbb{A}_1(u_1) = \left(A^i_{\alpha j} \,u^j_1\right),
\qquad \mathbb{A}_l(u_1,...,u_l) = \left(
\begin{array}{c}
\mathbb{A}_{l-1}(u_1,...,u_{l-1})
\\
A^i_{\alpha j} \,u^j_l
\end{array}
\right),
\]
see \cite[\S 5]{Cartan1}, \cite[Def. 11.4]{Olver95} for the full discussion. 
The system of forms $\omega^k$ is {\it involutive}  if (\ref{CT}) is an equality, \cite[\S 6]{Cartan1}, \cite[Def. 11.7]{Olver95}.

Cartan's fundamental theorems, 
\cite[\S\S 16, 22--24]{Cartan1}, \cite{Cartan4}, 
\cite[\S\S 16, 19, 20, 25,26]{Vasil'eva},
\cite[\S\S 14.1--14.3]{Stormark}, 
state that for a Lie pseudo-group there exists a set of {\sc mc} forms whose structure equations satisfy the compatibility and involutivity conditions; conversely, if Eqs. (\ref{SE}), (\ref{dUs}) meet the compatibility conditions (\ref{compatibility_conditions_SE}), (\ref{compatibility_conditions_dUs}) and the involutivity condition, then there exists a collection of 1-forms $\omega^1$, ... , $\omega^m$ and functions $U^1$, ... , $U^l$ which satisfy (\ref{SE}) and (\ref{dUs}).  Eqs. (\ref{Phi_omega}) then 
define local diffeomorphisms from a Lie pseudo-group.

\vskip 5 pt
\noindent
{\bf Example 1.}
Consider the bundle $J^2(\pi)$ of jets of the second order of the bundle $\pi$. A differential 1-form $\vartheta$ on $J^2(\pi)$ is called a {\it contact form} if it is annihilated by all 2-jets of local sections: $j_2(f)^{*}\vartheta = 0$. In the local coordinates every contact 1-form is a linear combination of the forms  $\vartheta_0 = du - u_{i}\,dx^i$, $\vartheta_i = du_i - u_{ij}\,dx^j$, $i, j \in \{1,...,n\}$, $u_{ji} = u_{ij}$. A local dif\-feo\-mor\-phism 
$\Delta : J^2(\pi) \rightarrow J^2(\pi)$, 
$\Delta : (x^i,u,u_i,u_{ij}) \mapsto (\widetilde{x}^i,\widetilde{u},\widetilde{u}_i,\widetilde{u}_{ij})$,
is called a {\it contact trans\-for\-ma\-tion} if for every contact 1-form $\widetilde{\vartheta}$ the form $\Delta^{*}\widetilde{\vartheta}$ is also contact. We denote by ${\rm{Cont}}(J^2(\pi))$  the pseudo-group  of contact transformations on $J^2(\pi)$.
Consider the following 1-forms 
\[
\Theta_0 = a\, \vartheta_0,
\quad
\Theta_i = g_i\,\Theta_0 + a\,B_i^k\,\vartheta_k,
\quad
\Xi^i =c^i\,\Theta_0+f^{ik}\,\Theta_k+b_k^i\,dx^k,
\]
\begin{equation}
\Theta_{ij} = a\,B^k_i\, B^l_j\,(du_{kl}-u_{klm}\,dx^m) + s_{ij}\,\Theta_0+w_{ij}^{k}\,\Theta_k,
\label{LCF}
\end{equation}
\noindent
defined on $J^2(\pi)\times\EuScript{H}$, where $\EuScript{H}$ is an open subset of $\mathbb{R}^{(2 n+1)(n+3)(n+1)/3}$ with local coordinates $(a$, $b^i_k$, $c^i$, $f^{ik}$, $g_i$, $s_{ij}$, $w^k_{ij}$, $u_{ijk})$, 
$i,j,k, \in \{1,...,n\}$,  such that $a\not =0$, $\det (b^i_k) \not = 0$, $f^{ik}=f^{ki}$, $s_{ij}=s_{ji}$, $w_{ij}^k=w_{ji}^k$,  $u_{ijk}=u_{ikj}=u_{jik}$, while $(B^i_k)$ is the inverse matrix for the matrix $(b^i_k)$.
As it is shown in \cite{Morozov2006}, the forms (\ref{LCF}) are {\sc mc} forms for ${\rm{Cont}}(J^2(\pi))$, that is, a local diffeomorphism
$\widehat{\Delta} : J^2(\pi) \times \EuScript{H} \rightarrow J^2(\pi) \times \EuScript{H}$
satisfies the conditions
$\widehat{\Delta}^{*}\, \widetilde{\Theta}_0 = \Theta_0$,
$\widehat{\Delta}^{*}\, \widetilde{\Theta}_i = \Theta_i$,
$\widehat{\Delta}^{*}\, \widetilde{\Xi}^i = \Xi^i$,
and $\widehat{\Delta}^{*}\, \widetilde{\Theta}_{ij} = \Theta_{ij}$
if and only if it is projectable on $J^2(\pi)$, and its projection
$\Delta : J^2(\pi) \rightarrow J^2(\pi)$ is a contact transformation.
The structure equations for ${\rm{Cont}}(J^2(\pi))$ have the form
\begin{eqnarray}
\fl\hspace{10pt}
d \Theta_0 &=& \Phi^0_0 \wedge \Theta_0 + \Xi^i \wedge \Theta_i,
\nonumber
\\
\fl\hspace{10pt}
d \Theta_i &=& \Phi^0_i \wedge \Theta_0 + \Phi^k_i \wedge \Theta_k,
\nonumber
\\
\fl\hspace{10pt}
d \Xi^i &=& \Phi^0_0 \wedge \Xi^i -\Phi^i_k \wedge \Xi^k
+\Psi^{i0} \wedge \Theta_0
+\Psi^{ik} \wedge \Theta_k,
\nonumber
\\
\fl\hspace{10pt}
d \Theta_{ij} &=& \Phi^k_i \wedge \Theta_{kj} 
+ \Phi^k_j \wedge \Theta_{ki}
- \Phi^0_0 \wedge \Theta_{ij}
+ \Upsilon^0_{ij} \wedge \Theta_0
+ \Upsilon^k_{ij} \wedge \Theta_k + \Xi^k \wedge \Theta_{ijk},
\nonumber
\end{eqnarray}
where the additional forms $\Phi^0_0$, $\Phi^0_i$, $\Phi^k_i$, $\Psi^{i0}$, $\Psi^{ij}$,
$\Upsilon^0_{ij}$, $\Upsilon^k_{ij}$, and $\Theta_{ijk}$ depend on differentials of the coordinates of $\EuScript{H}$.

\vskip 5 pt
 
\noindent
{\bf Example 2.}
Suppose $\EuScript{E}$ is a second-order differential equation in one dependent and $n$ independent variables. We consider $\EuScript{E}$ as a submanifold in $J^2(\pi)$.
Let ${\rm{Cont}}(\EuScript{E})$ be the group of contact symmetries for $\EuScript{E}$. It consists of all the contact transformations on $J^2(\pi)$ mapping $\EuScript{E}$ to itself. Let 
$\iota_0 : \EuScript{E} \rightarrow J^2(\pi)$ be an embedding, and 
$\iota = \iota_0 \times \rm{id} : \EuScript{E}\times \EuScript{H} \rightarrow J^2(\pi)\times \EuScript{H}$. The {\sc mc} forms of ${\rm{Cont}}(\EuScript{E})$ can be computed from  the forms  $\theta_0 = \iota^{*} \Theta_0$,  $\theta_i= \iota^{*}\Theta_i$, $\xi^i = \iota^{*}\Xi^i$, and $\theta_{ij}=\iota^{*}\Theta_{ij}$ by means of Cartan's method of equivalence, \cite{Cartan1,Cartan2,Cartan3,Cartan4,Gardner,Kamran,Olver95}, see details and examples in \cite{FelsOlver,Morozov2002,Morozov2006}.

\section{Cartan's structure of the contact symmetry pseudo-group for r-mdKP}
The r-th mdKP 
\[\fl
u_{tx} = -\frac{(3-r)\,(1-r)}{2}\,u_x^2 u_{xx} +\frac{r\,(3-r)}{2-r}\,u_x u_{xy} 
+\frac{3-r}{(2-r)^2}\,u_{yy}+\frac{(3-r)\,(1-r)}{2-r}\,u_y u_{xx},
\]
$r \in \mathbb{Z}\backslash \{2\}$,  was derived in \cite{Blaszak}. For a convenience of computations we use the following change of variables: 
\[
\tilde{t} = (3-r)\,t, 
\quad 
\tilde{x} = x,
\quad  
\tilde{y} = (2-r)\,y, 
\quad
\tilde{u} = -(1-r)\,u,
\]
where $r\not \in\{1, 2, 3\}$.  Then we have
\[
\tilde{u}_{\tilde{y}\tilde{y}} = \tilde{u}_{\tilde{t}\tilde{x}} +
\left(\frac{1}{2\,(1-r)}\,\tilde{u}_{\tilde{x}}^2+\tilde{u}_{\tilde{y}}\right)\,
\tilde{u}_{\tilde{x}\tilde{x}}
+\frac{r}{1-r}\,\tilde{u}_{\tilde{x}}\,\tilde{u}_{\tilde{x}\tilde{y}}.
\]
We drop tildes and denote $\kappa = \frac{r}{1-r}$; this yields %Eq. (\ref{main}).
\begin{equation}
u_{yy} = u_{tx}+\left(\textfrac{1}{2}\,(\kappa+1)\,u_x^2+u_y\right)\,u_{xx} +\kappa\,u_x\,u_{xy}.
\label{main}
\end{equation}
The exceptional cases $r=2$ and $r=3$ correspond to $\kappa = -2$ and $\kappa = -\frac{3}{2}$, respectively. We will not consider the case of $r=1$. The case of $\kappa = -1$ is exceptional, too, since $r \rightarrow \infty$ when $\kappa \rightarrow -1$. In the cases of $\kappa = 0$, 
$\kappa =1$, and $\kappa=-1$ Eq. (\ref{main}) gets the forms of the mdKP equation, \cite{Kuzmina,Krichever,Kupershmidt},  
%\begin{equation}
%\[
%u_{yy} = u_{tx}+\left(\textfrac{1}{2}\,u_x^2+u_y\right)\,u_{xx},
%\]
%\label{KZ}
%\end{equation}
the dBKP equation, \cite{Takasaki,KonopelchenkoAlonso},
%\begin{equation}
%\[
%u_{yy} = u_{tx}+\left(u_x^2+u_y\right)\,u_{xx} +u_x\,u_{xy}.
%\]
%\label{dBKP}
%\end{equation}
and the equation describing Lorentzian hyper-CR Einstein--Weil structures, \cite{Dunajski,FerapontovKhusnutdinova}.
%,
%\begin{equation}
%\[
%u_{yy} = u_{tx}+u_y\,u_{xx} -u_x\,u_{xy}.
%\]
%\label{Dunajski}
%\end{equation}

We use the method outlined in the previous section to compute {\sc mc} forms and structure equations of the pseudo-group of contact symmetries for Eq. (\ref{main}). In the next section, the analysis of integrable extensions will depend crucially on invariants of the pseudo-group. So in the computations we shall not suppress the invariants, as it was done in \cite{Morozov2008b}. We write out here the structure equations for the forms $\theta_0$, $\theta_j$, $\xi^j$ with $j \in \{ 1,2,3 \}$ only; the full sets of the structure equations are given in Appendix.

When $\kappa \not = -1$, we have 
%--------------------------------------------------------------------------------------
%     General case 
%--------------------------------------------------------------------------------------
\begin{eqnarray}
\fl 
d\theta_0 
&=& 
\eta_1 \wedge \theta_0+\xi^1 \wedge \theta_1
+\xi^2 \wedge \theta_2+\xi^3 \wedge \theta_3,
\nonumber
\\
\fl 
d\theta_1 
&=& 
\textfrac{1}{16}\,\left(
24\,\eta_1 
-24\,\theta_{22}
-12\,V\, \xi^1 
-12\,\xi^2 
-(12\,U-\kappa+4)\,\xi^3 
\right)
\wedge \theta_1
+\xi^1 \wedge \theta_{11}
+\xi^2 \wedge \theta_{12}
\nonumber
\\
\fl 
&+&
\left(
\textfrac{1}{8}\,(2\,\kappa^2+15\,\kappa+4)\,\theta_2 
+\kappa\,\theta_{23}
+(\kappa\,U+1-\textfrac{1}{8}\,\kappa\,(\kappa-10)) \,\xi^2
\right.
\nonumber
\\
\fl 
&+&
\left.
\textfrac{1}{16}\,(7\,\kappa+4)\,V\,\xi^3
\right) 
\wedge \theta_0
+
((\kappa+1)\,\theta_2+(\kappa+2)\,\xi^2) \wedge \theta_3
+\xi^3 \wedge \theta_{13},
\nonumber
\\
\fl 
d\theta_2
&=&
\textfrac{1}{2}\,\left(
\eta_1 
-\theta_{22}
-\textfrac{1}{2}\,V\,\xi^1
-\xi^2 
-\left(U -\textfrac{7}{8}\,\kappa-\textfrac{1}{2}\right)\,\xi^3 
\right)
\wedge \theta_2
+\xi^1 \wedge \theta_{12} 
+\xi^2 \wedge \theta_{22}
+\xi^3 \wedge \theta_{23},
\nonumber
\\
\fl 
d\theta_3
&=&
\left(\eta_1 -\theta_{22}-\textfrac{1}{2}\,V\,\xi^1-\xi^2
-\left(U-\textfrac{1}{4}\,\kappa\right)\,\xi^3 
\right)
\wedge \theta_3
+\textfrac{1}{2}\,(\kappa+2)\,\xi^2 \wedge \theta_2 
+\xi^1 \wedge \theta_{13}
\nonumber
\\
\fl 
&+&
\textfrac{1}{8}\,(\kappa-4)\,
\left(\theta_{22}+\xi^2+\textfrac{1}{8}\,(8\,U-\kappa-4)\,\xi^3\right)\wedge\theta_0
+\xi^2 \wedge \theta_{23}
+\xi^3 \wedge \theta_{12},
\nonumber
\\
\fl 
d\xi^1
&=&
-\textfrac{1}{2}\,\left(
\eta_1-3\,\theta_{22}+3\,\xi^2-\left(3\,U-\textfrac{1}{8}\,(\kappa-4)\right)\,\xi^3 
\right)
\wedge \xi^1,
\nonumber
\\
\fl 
d\xi^2
&=&
\textfrac{1}{2}\,\left(
\eta_1+\theta_{22}+\textfrac{1}{2}\,V\,\xi^1
+\left(U+\textfrac{1}{8}\,\kappa+\textfrac{3}{2}\right)\,\xi^3 
\right)
\wedge \xi^2
+\left(\textfrac{1}{8}\,(\kappa-4)\,\theta_0 -\theta_3 \right) \wedge \xi^1
-\theta_2 \wedge \xi^3,
\nonumber
\\
\fl 
d\xi^3
&=&
\left(
\theta_{22} +\textfrac{1}{2}\,V\,\xi^1 +\xi^2 
\right) 
\wedge \xi^3
-(\kappa+2)\,\left(\theta_2 +\xi^2 \right) \wedge \xi^1,
\nonumber
\end{eqnarray}
where the invariants
\begin{eqnarray}
\fl 
U &=& \frac{u_{xxy}+u_x\,u_{xxx} }{u_{xx}^{2}}, 
\label{U_general}
\\
\fl 
V &=& \frac{
u_{xxx}\left(
2\,u_{txx}
+((\kappa+3)\,u_x^2+2\,u_y)\,u_{xxx}
+2\,(\kappa+2)\,(u_x\,(u_{xxy}+u_{xx}^2)+u_{xx}\,u_{xy})
\right)}{u_{xx}^{4}}, 
\label{V_general}
\end{eqnarray}
satisfy
\begin{eqnarray}
\fl
d U &=&
\theta_2-U\,\theta_{22}-\eta_2
-\textfrac{1}{64}\,V\left(64\,U+7\,\kappa^2+38\,\kappa-8\right)\,\xi^1
-\textfrac{1}{16}\,\left(8\,U-\kappa-6\right)\,\xi^2
\nonumber
\\
\fl
&-&
\left(2\,U^2 -(2\,\kappa+1)\,U-\textfrac{3}{4}\,V -\kappa-1\right)\,\xi^3,
\nonumber
\\
\fl
d V 
&=& 
-\textfrac{1}{4}\,(\kappa-4)\,\theta_0
+\textfrac{1}{4}\,\left(8\,(\kappa+2)\,U+3\,(\kappa^2+6\,\kappa+8)\right)\,\theta_2
+2\,\theta_3
-\textfrac{3}{2}\,V\,\theta_{22}
%\label{d_invariants_general}
\nonumber
\\
\fl
&+&
2\,(\kappa+2)\,\theta_{23}
-2\,\eta_3
+\textfrac{1}{2}\,V\,\eta_1
+\left(2\,(\kappa+2)\,(2\,U+1)-V\right)\,\xi^2
\nonumber
\\
\fl
&-&
\textfrac{1}{32}\,V\left(80\,U+7\,\kappa^2+4\,\kappa-64\right)\,\xi^3,
\nonumber  
\end{eqnarray}
and where 
\begin{eqnarray}
\fl
\theta_0
&=&
\frac{u_{xxx}^2}{u_{xx}^3}\,\vartheta_0,
\nonumber
\\
\fl
\theta_1
&=&
\frac{u_{xxx}^3}{{u_{xx}}^6}\,\left(
\kappa\,(u_x\,u_{xx}+u_{xy})\,\vartheta_{0}
+\vartheta_1
+\left(\textfrac{1}{2}\,(\kappa+3)\,{u_x}^2+u_y\right)\,\vartheta_2
+(\kappa+2)\,u_x\,\vartheta_3
\right),
\nonumber
\\
\fl
\theta_2
&=&
\frac{u_{xxx}}{u_{xx}^2}\,\vartheta_2,
\nonumber
\\
\fl
\theta_3
&=&
\frac{u_{xxx}^2}{u_{xx}^4}\,
\left(
(\kappa-4)\,u_{xx}\,\vartheta_{0}
+u_x\,\vartheta_2
+\vartheta_3
\right),
\nonumber
\\
\fl
\theta_{22}
&=&
\frac{1}{u_{xx}}\,\vartheta_{22},
\label{MC_forms_general}
\\
\fl
\xi^1
&=&
\frac{u_{xx}^3}{u_{xxx}}\,dt,
\nonumber
\\
\fl
\xi^2
&=&
\frac{u_{xxx}}{u_{xx}}\,\left(
\left(\textfrac{1}{2}\,(\kappa+1)\,u_{xx}^2-u_y\right)\,dt
+dx-u_x\,dy
\right),
\nonumber
\\
\fl
\xi^3
&=&
-u_{xx}\,u_x\,(\kappa+2)\,dt+u_{xx}\,dy,
\nonumber
\\
\fl
\eta_1 
&=&
2\,\frac{d u_{xxx}}{u_{xxx}}
-3\,(\theta_{22}-\xi^2)
-\frac{1}{2\,u_{xx}^3}\,\left(
3\,V\,u_{xx}^3
-4\,(\kappa+3)\,(u_x\,u_{xx}+u_{xy})\,u_{xxx}
\right)\,\xi^1
\nonumber
\\
\fl
&-&\left(3\,U+\textfrac{1}{8}(\kappa-4)\right)\,\xi^3,
\nonumber
\end{eqnarray}
with 
\begin{eqnarray}
\fl
\vartheta_0
&=&
du-u_t\,dt-u_x\,dx-u_y\,dy,
\nonumber
\\
\fl 
\vartheta_1
&=&du_t -u_{tt}\,dt-u_{tx}\,dx-u_{ty}\,dy,
\nonumber
\\
\fl
\vartheta_2
&=&
du_x -u_{tx}\,dt-u_{xx}\,dx-u_{xy}\,dy,
\nonumber
\\
\fl
\vartheta_3
&=&
du_y-u_{ty}\,dt-u_{xy}\,dx
-(u_{tx}+\left(\textfrac{1}{2}(\kappa+1)\,u_x^2+u_y\right)\,u_{xx} +\kappa\,u_x\,u_{xy})\,dy,
\nonumber
\\
\fl
\vartheta_{22}&=&du_{xx} -u_{txx}\,dt-u_{xxx}\,dx-u_{xxy}\,dy.
\nonumber
\end{eqnarray}
We need not explicit expressions for the other {\sc mc} forms in the sequel.

\vskip 5 pt

For $\kappa = -1$, we get
%--------------------------------------------------------------------------------------
%     kappa = - 1
%--------------------------------------------------------------------------------------
\begin{eqnarray}
\fl
d\theta_0
&=&
\eta_1 \wedge \theta_0+\xi^1 \wedge \theta_1+\xi^2 \wedge \theta_2
+\xi^3 \wedge \theta_3,
\nonumber
\\
\fl
d\theta_1
&=&
\textfrac{3}{2}\,\eta_1 \wedge \theta_1
+2\,\eta_2 \wedge \theta_3
+\eta_3 \wedge \theta_0
-\textfrac{1}{8}\,\left(
12\,\theta_{22}+16\,U\,\xi^1+12\,\xi^2+7\,\xi^3
\right) \wedge \theta_1
+\xi^1 \wedge \theta_{11}
\nonumber
\\
\fl
&+&
\xi^2 \wedge \theta_{12}
+\xi^3 \wedge \theta_{13},
\nonumber
\\
\fl
d\theta_2
 &=&
\textfrac{1}{8}\,\left(
4\,(\eta_1-\theta_{22} -\xi^2)-3\,\xi^3
\right)\wedge \theta_2
+\xi^1 \wedge \theta_{12}+\xi^2 \wedge \theta_{22}+\xi^3 \wedge \theta_{23},
\nonumber
\\
\fl
d\theta_3
&=&
\left(
\eta_1-\theta_{22}-U\,\xi^1-\xi^2-\textfrac{5}{8}\,\xi^3 
\right)\wedge \theta_3
+\eta_2 \wedge \theta_2
-\textfrac{5}{8}\,(\theta_{22}+\xi^2) \wedge \theta_0 
+\xi^1 \wedge \theta_{13}
\nonumber
\\
\fl
&+&
\xi^2 \wedge \theta_{23}
+\xi^3 \wedge \theta_{12},
\nonumber
\\
\fl
d\xi^1
&=&
\textfrac{1}{8}\,\left(
12\,\theta_{22}-4\,\eta_1+12\,\xi^2+7\,\xi^3
\right) \wedge \xi^1,
\nonumber
\\
\fl
d\xi^2
&=&
-\textfrac{1}{8}\,\left(
5\,\theta_0 + 8\,\theta_3\right) \wedge \xi^1
+\textfrac{1}{8}\,\left(
4\,(\eta_1 +\theta_{22})+3\,\xi^3
\right) \wedge \xi^2
-\left(\eta_2 +\theta_2 \right)\wedge \xi^3,
\nonumber
\\
\fl
d\xi^3
&=&
-(2\,\eta_2 +\theta_2) \wedge \xi^1
+(\theta_{22} 
+U\,\xi^1+\xi^2) \wedge \xi^3,
\nonumber
\end{eqnarray}
where the invariant
\begin{equation}
\fl
U=\frac{
u_{xxx}\,(u_x\,u_{xxy}-u_{txx}-u_y\,u_{xxx}-u_{xy}\,u_{xx})
+u_{xxy}\,(u_{xxy}+u_{xx}^2)
}{u_{xx}^4}
-\frac{33}{64}
\label{U_exceptional}
\end{equation}
satisfies
\begin{equation}
\fl
d U=
-\textfrac{5}{8}\,\theta_0
+\textfrac{9}{8}\,\theta_2
-\textfrac{1}{8}\,(4\,U-7)\,\xi^2 
-\textfrac{3}{2}\,U\,(\theta_{22}+\xi^3)
+\theta_{23}
+\eta_3-\textfrac{1}{4}\,\eta_2
+\textfrac{1}{2}\,U\,\eta_1,
%\label{d_invariant_exceptional}
\nonumber
\end{equation}
and where
\begin{eqnarray}
\fl
\theta_0&=&\frac{u_{xxx}^2}{{u_{xx}}^3}\,\vartheta_{0},
\nonumber
\\
\fl
\theta_1&=&
\frac{u_{xxx}^3}{u_{xx}^6}\,\vartheta_1
+\frac{u_{xxx}}{64\,u_{xx}^6}\,
\left(
64\,(u_y\,u_{xxx}^2+u_{xxy}^2)
+u_{xx}^2\,(9\,u_{xx}^2-48\,u_{xxy})
\right)\,\vartheta_2
\nonumber
\\
\fl
&-&
\frac{u_{xxx}^2}{u_{xx}^7}\left(
u_{xxx}\,(u_{txx}+u_y\,u_{xxx}-u_x\,u_{xxy}
+2\,u_{xx}\,u_{xy})
-u_{xxy}\,(u_{xxy}+2\,\,u_{xx}^2)
+\textfrac{57}{64}\,u_{xx}^4
\right)\,\vartheta_{0}
\nonumber
\\
\fl
&-&
\frac{u_{xxx}^2}{4\,u_{xx}^6}\,\left(
8\,u_{xxy}+4\,u_x\,u_{xxx}
-3\,u_{xx}^2
\right)\,\vartheta_3,
\nonumber
\\
\fl
\theta_2
&=&
\frac{u_{xxx}}{u_{xx}^2}\,\vartheta_2,
\nonumber
\\
\fl
\theta_3
&=&
\frac{u_{xxx}}{8 u_{xx}^4}\,\left(
8\,u_{xxx}\,\vartheta_3
-5\,u_{xx}\,u_{xxx}\,\vartheta_{0}
-(8\,u_{xxy}-3\,u_{xx}^2)\vartheta_2
\right),
\nonumber
\\
\fl
\theta_{22}
&=&
\frac{1}{u_{xx}}\,\vartheta_{22},
\label{MC_forms_exceptional}
\\
\fl
\xi^1
&=&
\frac{u_{xx}^3}{u_{xxx}}\,dt,
\nonumber
\\
\fl
\xi^2
&=&
\frac{1}{u_{xx} u_{xxx}}
\left(
u_{xxy}\,\left(u_{xxy}-\textfrac{3}{4}\,u_{xx}^2\right)
+u_{xxx}\,(u_x\,u_{xxy}-u_y\,u_{xxx})
\right)\,dt
\nonumber
\\
\fl
&-&
\textfrac{3}{64}\,u_{xx}^2\,\left(8\,u_x\,u_{xxx}-3\,u_{xx}^2\right)
+
\frac{u_{xxx}}{u_{xx}}\,dx
+\frac{(8\,u_{xxy}-3\,u_{xx}^2)}{8\,u_{xx}}\,dy,
\nonumber
\\
\fl
\xi^3
&=&
\frac{u_{xx}}{4 u_{xxx}}\, \left(
8\,u_{xxy}
+4\,u_x\,u_{xxx}
-3\,u_{xx}^2
\right)\,dt
+u_{xx}\,dy,
\nonumber
\\
\fl
\eta_1
&=&
2\,\frac{d u_{xxx}}{u_{xxx}}
-3\,(\theta_{22}+\xi^2)
+\frac{1}{2 u_{xx}^3}\,\left(
(4\,U+3)\,u_{xx}^3\,
+8\,(u_{xy}\,u_{xxx}-u_{xx}\,u_{xxy})
\right)\,\xi^1
-\textfrac{1}{2}\,\xi^3.
\nonumber
\end{eqnarray}

\section{Integrable extensions}

In \cite[\S 6]{BryantGriffiths}, the definition of integrable extension of an exterior differential system is designed to study finite-dimensional coverings. In general case, coverings of {\sc pde}s with three or more independent variables are infinite-dimensional, \cite{Marvan1992}. To cope with infinite-dimensional coverings we use a  natural generalization of the definition. 

Suppose $\mathfrak{G}$ is a Lie pseudo-group on a manifold $M$ and $\omega^1$, ... , $\omega^m$ are its {\sc mc} forms with structure equations (\ref{SE}), (\ref{dUs}). Consider a system of equations
\begin{eqnarray}
\fl
d\tau^q &=& 
D^q_{\rho r} \, \eta^\rho \wedge \tau^r 
+ 
E^q_{r s} \, \tau^r \wedge \tau^s 
+
F^q_{r \beta} \, \tau^r \wedge \pi^\beta
+
G^q_{r j} \, \tau^r \wedge \omega^j
+
H^q_{\beta j} \, \pi^\beta \wedge \omega^j
\nonumber
\\
\fl
&+&
I^q_{j k} \, \omega^j \wedge \omega^k,
\label{extra_SE} 
\\
\fl
d V^\epsilon &=&  J^\epsilon_j \,\omega^j 
+  K^\epsilon_q \, \tau^q,
\label{dVs}
\end{eqnarray}
with unknown 1-forms $\tau^q$, $q \in \{1,...,Q\}$, $\eta^\rho$, $\rho \in \{1,...,R\}$,
and unknown functions $V^\epsilon$, $\epsilon \in \{1,...,S\}$ for some $Q, R, S \in \mathbb{N}$. The coefficients $D^q_{\rho r}$, ..., $K^\epsilon_q$ in (\ref{extra_SE}), (\ref{dVs}) are supposed to be functions of $U^q$ and $V^\gamma$.

\vskip 5 pt

\noindent
{\bf Definition 1.}
System (\ref{extra_SE}), (\ref{dVs}) is an {\it integrable extension} of system (\ref{SE}), (\ref{dUs}), if Eqs. (\ref{extra_SE}), (\ref{dVs}), (\ref{SE}), and (\ref{dUs}) together satisfy the compatibility and involutivity conditions. 

\vskip 5 pt

In this case from Cartan's third fundamental theorem for Lie pseudo-groups it follows that there exists a set of forms $\tau^q$ and functions $V^\epsilon$ which are solutions to Eqs. (\ref{extra_SE}) and (\ref{dVs}). Then $\tau^q$, $V^\epsilon$ together with $\omega^i$, $U^q$ define a Lie pseudo-group %$\mathfrak{H}$ 
on a manifold $N \cong M \times \mathbb{R}^Q$.

\vskip 5 pt

\noindent
{\bf Definition 2.}
The integrable extension is called {\it trivial}, if there is a change of variables on $N$ such that in the new variables the coefficients $F^q_{r \beta}$, $G^q_{r j}$, $H^q_{\beta j}$, $I^q_{j k}$, and $J^\epsilon_j$ are equal to zero, while the coefficients $D^q_{\rho r}$, $E^q_{r s}$, and $K^\epsilon_q$ are independent of $U^q$. Otherwise, the integrable extension is called {\it non-trivial}.

\vskip 5pt 

Let $\theta^\alpha_I$ and $\xi^j$ be a set of {\sc mc} forms of the symmetry pseudo-group $\mathrm{Cont}(\EuScript{E})$ of a {\sc pde} $\EuScript{E}$ such that $\xi^1 \wedge ... \wedge \xi^n \not = 0$ on any solution manifold of 
$\EuScript{E}$, while $\theta_I^\alpha$ are contact forms. 
We take the following reformulation of the  definition (\ref{WE_def}) of a covering. 

\vskip 5 pt

\noindent
{\bf Definition 3.}  A non-trivial integrable extension of the form
\begin{equation}
d \omega^q =\Pi^q_r \wedge \omega^r + \xi^j \wedge \Omega^q_j 
\label{contact_ie}
\end{equation}
is referred to as a {\it contact integrable extension} ({\sc cie}) of the structure equations of $\mathrm{Cont}(\EuScript{E})$ if

(1)\,\,\, $\Pi^q_r$ are some non-trivial differential 1-forms,

(2)\,\,\, $\Omega^q_j \equiv 0 \,\,\,{\mathrm{mod}}\,\,\,\theta^{\alpha}_I, \omega^q_j$ 
for some additional 1-forms $\omega^q_j$.

\vskip 7 pt

Since (\ref{contact_ie}) is integrable extension, Cartan's theorem yields existence of the forms $\omega^q$ satisfying (\ref{contact_ie}). From \cite[Ch. IV, Prop. 5.10]{BCGGG} it follows that the forms $\omega^q$ define a system of {\sc pde}s. This system is a covering for $\EuScript{E}$.

We apply this construction to the structure equations (\ref{SE_general}) and (\ref{SE_exceptional}) of the symmetry pseudo-group of the r-th mdKP equation.
We restrict our analysis to {\sc cie}s of the form 
\begin{eqnarray}
d \omega_0 
&=& 
\left(
\sum \limits_{i=0}^3 A_i \,\theta_i 
+ %\sum \limits_{1\le i \le j \le 3,\,\, (i,j)\not = (3,3)} 
  \sum {}^{*} B_{ij}\,\theta_{ij}
+ \sum \limits_{s=1}^7 C_s\,\eta_s 
+ \sum \limits_{j=1}^3 D_j\,\xi^j 
+ E\,\omega_1  
\right) \wedge \omega_0 
\nonumber
\\
&+& 
\sum \limits_{j=1}^3 \left(
\sum \limits_{k=0}^3 F_{jk}\,\theta_k + G_j\,\omega_1
\right) \wedge \xi^j,
\label{ie_main}
\end{eqnarray}
where $\sum {}^{*}$ means summation for all $i,j \in \mathbb{N}$ such that $1\le i \le j \le 3$, $(i,j)\not = (3,3)$, and consider two cases:

1) the coefficients $A_i$ to $G_j$ in (\ref{ie_main}) depend on the invariants of the symmetry pseudo-group of Eq. (\ref{main}) only, i.e., on the $U$ and $V$ when $\kappa\not = -1$ and on $U$ when $\kappa =-1$;

2) the coefficients $A_i$ to $G_j$ depend also on one additional invariant, say $W$. In this case, the differential of this new invariant satisfies the following equation
\begin{equation}
d W = \sum \limits_{i=0}^3 H_i \,\theta_i 
+ \sum {}^{*} %\limits_{1\le i \le j \le 3,\,\, (i,j)\not = (3,3)} 
    I_{ij}\,\theta_{ij}
+ \sum \limits_{s=1}^7 J_s\,\eta_s 
+ \sum \limits_{j=1}^3 K_j\,\xi^j 
+ \sum \limits_{q=0}^1 L_q\,\omega_q,
\label{dW}
\end{equation}
where the coefficients $H_i$ to $L_q$ are functions of $U$, $V$, $W$ when $\kappa\not = -1$ or of $U$ and $W$ when $\kappa= -1$.

The requirements of Defintions 1 and 3 yield over-determined systems for the coefficients $A_i$ to $G_j$ in the first case and $A_i$ to $L_q$ in the second case. The results of analysis of these systems are summarized in the following theorems.

%=======================================================================================
\vskip 5 pt
\noindent
{\bf Theorem 1.} 
{\it 
Let $\kappa \not = -1$ and suppose that the coefficients $A_i$ to $G_j$ in (\ref{ie_main})  depend on the invariants $U$, $V$ defined by (\ref{U_general}) and (\ref{V_general}). Then 
{\sc cie} (\ref{ie_main}) of (\ref{SE_general})  is contact-equivalent to the following equation:
}
\begin{eqnarray}
\fl
d\omega_0
&=&
\left(
\omega_1
+\textfrac{1}{2}\,(\eta_1+\theta_{22})
+\textfrac{1}{4}\,V\,\xi^1
+\textfrac{1}{16}\,(8\,U+\kappa+12)\,\xi^3 
\right) \wedge \omega_0
+\omega_1 \wedge \xi^2
+\theta_2 \wedge \xi^3
\nonumber
\\
\fl
&-&\left(
\textfrac{1}{8}\,(\kappa-4)\,\theta_0 
-\theta_3
\right)\wedge \xi^1.
\label{extra_se_1}
\end{eqnarray}

%=======================================================================================

\vskip 5 pt
\noindent
{\bf Theorem 2.} 
{\it 
For $\kappa \not \in \{-3, -1\}$ suppose that the coefficients $A_i$ to $G_j$ in (\ref{ie_main}) and $H_i$ to $L_q$ in (\ref{dW}) depend on the invariants $U$, $V$, and $W$, where $U$ and $V$ are  defined by (\ref{U_general}), (\ref{V_general}). Then  {\sc cie} (\ref{ie_main}), (\ref{dW}) of (\ref{SE_general}) is contact-equivalent to the following equations:
}
\begin{eqnarray}
\fl
d\omega_0
&=&
\left(
\omega_1
+\textfrac{1}{2}\,(\eta_1+\theta_{22})
+\textfrac{1}{4}\,\left(V-2\,W\,(W-2)\right)\xi^1
+\textfrac{1}{16}\,(8\,(U-W)+\kappa+12)\,\xi^3
\right) \wedge \omega_0
\nonumber
\\
\fl
&+&
\left(
W^2\,\omega_1
-\textfrac{1}{8}\,(\kappa-4)\,\theta_0
+W\,\theta_2
+\theta_3
\right) \wedge \xi^1
+\omega_1 \wedge \xi^2
+\left(
W\,\omega_1
+\theta_2
\right) \wedge \xi^3,
\label{extra_se_2}
\\
\fl
dW
&=&\textfrac{1}{16}\,W\,\left(
8\,(\eta_1-\theta_{22}-\kappa\,\omega_1-(\kappa+1)\,\omega_2)
-4\,\left(V-2\,W\,(W-2\,(\kappa+1))\right)\,\xi^1
\right.
\nonumber
\\
\fl
&-&
\left.
(8\,(U-W)+5\,(3\kappa+4))\,\xi^3
\right).
\label{dW_2}
\end{eqnarray}
{\it
In the case of  $\kappa=-3$, every {\sc cie} (\ref{ie_main}), (\ref{dW}) of (\ref{SE_general})  is contact-equivalent either to (\ref{extra_se_2}), (\ref{dW_2}) or to the following equations: 
}
\begin{eqnarray} 
\fl
d\omega_0
&=&
\left(\omega_1+\textfrac{1}{2}\,(\eta_1+\theta_{22})
+\textfrac{1}{4}\,V\,\xi^1 
+\textfrac{3}{2}\,\xi^2
+\textfrac{1}{16}\,(8\,(U+2\,W)+9)\,\xi^3
\right) \wedge \omega_0
\nonumber
\\
\fl
&+&\left(\left(W+\textfrac{7}{8}\right)\,\theta_0+\theta_3\right) \wedge \xi^1
+\omega_1 \wedge \xi^2
+\theta_2 \wedge \xi^3,
\label{extra_se_3}
\\
\fl
dW
&=&
W\,\left(
\omega_1-\theta_{22}-\textfrac{1}{2}\,V\,W\,\xi^1
-(U-W-1)\,\xi^3
\right).
\label{dW_3}
\end{eqnarray}

%=======================================================================================

\vskip 5 pt
 
\noindent
{\bf Theorem 3.}
{\it 
In the case of $\kappa=-1$ there are no {\sc cie}s of (\ref{SE_exceptional})  with the coefficients $A_i$ to $G_j$ depending on the invariant $U$ from (\ref{U_exceptional}) only. Every {\sc cie} (\ref{ie_main}), (\ref{dW}) of (\ref{SE_exceptional}) whose coefficients $A_i$ to $L_q$ depend on $U$ and $W$ is contact-equivalent to the following equations: 
}
\begin{eqnarray}
\fl
d\omega_0
&=&
\left(\omega_1
+\textfrac{1}{2}\,(\eta_1+\theta_{22}-W\,(W-2)\,\xi^1)
-\textfrac{1}{8}\,(4\,W-7)\,\xi^3
\right) \wedge \omega_0
\nonumber
\\
\fl
&+&
\left(W^2\,\omega_1
+\textfrac{5}{8}\,\theta_0
+W\,\theta_2
+\theta_3\right) \wedge \xi^1 
+\omega_1 \wedge \xi^2
+(W\,\omega_1 +\theta_2) \wedge \xi^3,
\label{extra_se_4}
\\
\fl
dW&=&
\eta_2+U\,W\,\xi^1
+\textfrac{1}{2}\,W\,\left(\eta_1-\theta_{22}
-(W+1)\,\xi^2
+(2\,U-W)\,\xi^3\right).
\label{dW_4}
\end{eqnarray}

%=======================================================================================
Since forms (\ref{MC_forms_general}) in (\ref{extra_se_1}) -- (\ref{dW_3}) and forms (\ref{MC_forms_exceptional}) in (\ref{extra_se_4}), (\ref{dW_4}) are known, it is easy to find the forms $\omega_0$ explicitly:

%=======================================================================================
\vskip 5 pt

{\bf Theorem 4.}
{\it 
We have the following solutions to Eqs. (\ref{extra_se_1}), (\ref{extra_se_2}), (\ref{extra_se_3}), and (\ref{extra_se_4}) up to a contact equivalence:

1) Eq. (\ref{extra_se_1})  yields
}
\begin{equation}
\fl
\omega_0 = \frac{u_{xxx}}{u_{xx}\,v_x}\,
\left(
d v-\left(\textfrac{1}{2}\,(\kappa+1)\,u_x^2-u_y\right)\,v_x\,dt-v_x\,dx
+u_x\,v_x\,dy
\right).
\label{WE_1}
\end{equation}

{\it 
2) Integrating Eqs. (\ref{extra_se_2}), (\ref{dW_2}) yields the following results:

when $\kappa \not \in \{-2, -\textfrac{3}{2}, -1\}$, we have 
}
\begin{eqnarray}
\fl
\omega_0 &=& \frac{u_{xxx}}{u_{xx}\,v_x}\,\left(
dv
-
\left(
(2\kappa+3)^{-1}\,v_x^{2\,(\kappa+1)}
-\lambda\,u_x\,v_x^{\kappa+1}
+\left(\textfrac{1}{2}\,(\kappa+1)\,u_x^2-\,u_y\right)
\right)\,v_x\,dt
\right.
\nonumber
\\
\fl
&-&
\left.
v_x\,d x
-\left(
(\kappa+2)^{-1}\,v_x^{\kappa+1}
-u_x
\right)\,v_x\,dy
\right),
\label{WE_2}
\end{eqnarray}

{\it 
when $\kappa = -2$, we have
}
\begin{eqnarray} 
\fl
\omega_0 
&=& 
\frac{u_{xxx}}{u_{xx} v_x}\,\left(
d v
+\left(
v_x^{-1}
+u_x
+\left(\textfrac{1}{2}\,u_x^2+u_y\right)\,v_x
\right)\,dt
-\left(\ln \vert v_x \vert - u_x\,v_x\right)\,dy
%\right.
%\nonumber
%\\
%&-&
%\left.
-
v_x\,dx
\right),
\label{WE_3} 
\end{eqnarray}

{\it 
and when $\kappa=-\textfrac{3}{2}$, we have 
}
\begin{eqnarray} 
\fl
\omega_0
&=&
\frac{u_{xxx}}{u_{xx}\,v_x}\,\left(
d v
+\left(\left(\textfrac{1}{4}\,u_x^2+u_y\right)\,v_x
+u_x\,\sqrt{\vert v_x \vert} - \ln \vert v_x \vert
\right)\,dt
+\left(u_x\,v_x - 2 \,\sqrt{\vert v_x \vert}\right)\,d y
\right.
\nonumber
\\
\fl
&-&
\left.
v_x\,dx
\right).
\label{WE_4} 
\end{eqnarray}

{\it 
3) From Eqs. (\ref{extra_se_3}), (\ref{dW_3}) we get}
\begin{eqnarray} 
\fl
\omega_0 &=& \frac{u_{xxx}}{u_{xx}\,v_x}\,
\left(
dv
+\left(u+(u_x^2+u_y)\,v_x\right)\,dt
-v_x\,dx
+(v_x\,u_x+x)\,dy
\right).
\label{WE_5} 
\end{eqnarray}

{\it 
4) Eqs. (\ref{extra_se_4}), (\ref{dW_4}) give}
\begin{eqnarray}
\fl
\omega_0 
&=& \frac{u_{xxx}}{u_{xx}\,v_x}\,\left(
dv
+v_x\left(u_y+\lambda\,u_x-\lambda^2\right)\,dt
-v_x\,dx
+v_x\,(u_x-\lambda)\,dy
\right),
\quad \lambda \in \mathbb{R}
\label{WE_6}
\end{eqnarray}

\vskip 5 pt
%=======================================================================================

Forms (\ref{WE_1}) -- (\ref{WE_6}) define infinite-dimensional coverings of Eq. (\ref{main}) with the fibre variables $v_j$, $j \in \mathbb{N} \cup \{0\}$, and the following extended total derivatives, respectively:
\begin{equation}
\fl
\left\{
\begin{array}{lll}
\widetilde{D}_t &=&\bar{D}_t + \sum \limits_{j=0}^{\infty}\,\widetilde{D}_x^j 
\left( 
\left(
\textfrac{1}{2}\,(\kappa+1)\,u_x^2-u_y
\right)\,v_1
\right)\,\frac{\partial}{\partial v_j},
\\
\widetilde{D}_x &=& \bar{D}_x + \sum \limits_{j=0}^{\infty} v_{j+1}\,\frac{\partial}{\partial v_j}, 
\\

\widetilde{D}_y &=&\bar{D}_y - \sum \limits_{j=0}^{\infty}\,\widetilde{D}_x^j 
\left(
u_x\,v_1
\right)\,\frac{\partial}{\partial v_j},
\end{array}
\right.
\label{cov_1}
\end{equation}
%-------------------------------------------------------------------------------------
\begin{equation}
\fl
\left\{
\begin{array}{lll}
\widetilde{D}_t &=&\bar{D}_t + \sum \limits_{j=0}^{\infty}\,\widetilde{D}_x^j 
\left( 
(2\kappa+3)^{-1}\,v_1^{2\,\kappa+3}
-u_x\,v_1^{\kappa+2}
+\left(\textfrac{1}{2}\,(\kappa+1)\,u_x^2-\,u_y\right)
\,v_1
\right)\,\frac{\partial}{\partial v_j},
\\
\widetilde{D}_x &=& \bar{D}_x + \sum \limits_{j=0}^{\infty} v_{j+1}\,\frac{\partial}{\partial v_j}, 
\\
\widetilde{D}_y &=&\bar{D}_y + \sum \limits_{j=0}^{\infty}\,\widetilde{D}_x^j 
\left(
(\kappa+2)^{-1}\,v_1^{\kappa+2}-u_x\,v_1
\right)
\frac{\partial}{\partial v_j},
\end{array}
\right.
\label{cov_2}
\end{equation}
%-------------------------------------------------------------------------------------
\begin{equation}
\fl
\left\{
\begin{array}{lll}
\widetilde{D}_t &=&\bar{D}_t - \sum \limits_{j=0}^{\infty}\,\widetilde{D}_x^j 
\left( 
v_1^{-1}  
+u_x
+\left(\textfrac{1}{2}\,u_x^2+u_y\right)\,v_1
\right)\,\frac{\partial}{\partial v_j},
\\ 
\widetilde{D}_x &=& \bar{D}_x + \sum \limits_{j=0}^{\infty} v_{j+1}\,\frac{\partial}{\partial v_j}, 
\\ 
\widetilde{D}_y &=&\bar{D}_y + \sum \limits_{j=0}^{\infty}\,\widetilde{D}_x^j 
\left(
\ln \vert v_1 \vert  - u_x\,v_1
\right)\,\frac{\partial}{\partial v_j},
\end{array}
\right.
\label{cov_3}
\end{equation}
%-------------------------------------------------------------------------------------
\begin{equation}
\fl
\left\{
\begin{array}{lll}
\widetilde{D}_t &=&\bar{D}_t - \sum \limits_{j=0}^{\infty}\,\widetilde{D}_x^j 
\left( 
\left(\textfrac{1}{4}\,u_x^2+u_y\right)\,v_1
+u_x\,\sqrt{\vert v_1 \vert} - \ln \vert v_1 \vert
\right)\,\frac{\partial}{\partial v_j},
\\ 
\widetilde{D}_x &=& \bar{D}_x + \sum \limits_{j=0}^{\infty} v_{j+1}\,\frac{\partial}{\partial v_j}, 
\\
\widetilde{D}_y &=&\bar{D}_y + \sum \limits_{j=0}^{\infty}\,\widetilde{D}_x^j 
\left(
2 \,\sqrt{\vert v_1 \vert }-u_x\,v_1
\right)
\frac{\partial}{\partial v_j},
\end{array}
\right.
\label{cov_4}
\end{equation}
%-------------------------------------------------------------------------------------
\begin{equation}
\fl
\left\{
\begin{array}{lll}
\widetilde{D}_t &=&\bar{D}_t - \sum \limits_{j=0}^{\infty}\,\widetilde{D}_x^j 
\left( 
u+(u_x^2+u_y)\,v_1
\right)\,\frac{\partial}{\partial v_j},
\nonumber
\\
\widetilde{D}_x &=& \bar{D}_x + \sum \limits_{j=0}^{\infty} v_{j+1}\,\frac{\partial}{\partial v_j}, 
\\
\widetilde{D}_y &=&\bar{D}_y - \sum \limits_{j=0}^{\infty}\,\widetilde{D}_x^j 
\left(
v_1\,u_x+x
\right)
\frac{\partial}{\partial v_j},
\end{array}
\right.
\label{cov_5}
\end{equation}
%-------------------------------------------------------------------------------------
\begin{equation}
\fl
\left\{
\begin{array}{lll}
\widetilde{D}_t &=&\bar{D}_t - \sum \limits_{j=0}^{\infty}\,\widetilde{D}_x^j 
\left( 
(u_y+\lambda\,u_x-\lambda^2)\,v_1
\right)\,\frac{\partial}{\partial v_j},
\\ 
\widetilde{D}_x &=& \bar{D}_x + \sum \limits_{j=0}^{\infty} v_{j+1}\,\frac{\partial}{\partial v_j}, 
\\
\widetilde{D}_y &=&\bar{D}_y - \sum \limits_{j=0}^{\infty}\,\widetilde{D}_x^j 
\left(
(u_x-\lambda)\,v_1
\right)
\frac{\partial}{\partial v_j}.
\end{array}
\right.
\label{cov_6}
\end{equation}

The covering (\ref{cov_2}) was found in the cases of $\kappa=0$ and $\kappa=1$ in \cite{ChangTu} and \cite{KonopelchenkoAlonso}, respectively. The covering (\ref{cov_6}) was found in \cite{Dunajski}. In \cite{Morozov2008b}, we construct the coverings (\ref{cov_1}), (\ref{cov_2}) for $\kappa \not \in \{-2, -\frac{3}{2}, -1\}$, and the covering (\ref{cov_6}) by means of another technique. Another covering of Eq. (\ref{main}) in the case of $\kappa = - 1$ was found in \cite{Morozov2008b}, too. As far as we know,  the coverings (\ref{cov_3}), (\ref{cov_4}), and (\ref{cov_5}) are new.

\section{Conclusion}

We have shown that the approach based on the contact integrable extensions allows one to derive coverings of {\sc pde}s with three independent variables from {\sc mc} forms of their symmetry pseudo-groups. We hope that this technique will be effective in studying a wide variety of {\sc pde}s with large symmetry pseudo-groups. Two important problems require further exploration. The first one is a connection between the restrictions like (\ref{ie_main}) for the form of a {\sc cie} and the number of the additional invariants with structure properties of the symmetry pseudo-groups. The second task is a refinement of the analysis of over-determined systems of {\sc pde}s for the coefficients of {\sc cie}s by means of Cartan's method. We hope to address these problems in future papers.

\section*{Appendix}
The structure equations of the symmetry pseudo-group of Eq. (\ref{main}) in the case of $\kappa \not = -1$ read
%----------------------------------------------------------------------------------
%       General case
%----------------------------------------------------------------------------------
{\small
\begin{eqnarray} 
\fl
d\theta_0 
&=& 
\eta_1 \wedge \theta_0+\xi^1 \wedge \theta_1
+\xi^2 \wedge \theta_2+\xi^3 \wedge \theta_3,
\nonumber
\\ \fl
d\theta_1 
&=& 
\textfrac{1}{16}\,\left(
24\,\eta_1 
-24\,\theta_{22}
-12\,V\, \xi^1 
-12\,\xi^2 
-(12\,U-\kappa+4)\,\xi^3 
\right)
\wedge \theta_1
+\xi^1 \wedge \theta_{11}
+\xi^2 \wedge \theta_{12} 
\nonumber
\\ \fl
&+&
\left(
\textfrac{1}{8}\,(2\,\kappa^2+15\,\kappa+4)\,\theta_2 
+\kappa\,\theta_{23}
+(\kappa\,U+1-\textfrac{1}{8}\,\kappa\,(\kappa-10)) \,\xi^2
\right.
\nonumber
\\ \fl
&+&
\left.
\textfrac{1}{16}\,(7\,\kappa+4)\,V\,\xi^3
\right) 
\wedge \theta_0
+
((\kappa+1)\,\theta_2+(\kappa+2)\,\xi^2) \wedge \theta_3
+\xi^3 \wedge \theta_{13},
\nonumber
\\ \fl
d\theta_2
&=&
\textfrac{1}{2}\,\left(
\eta_1 
-\theta_{22}
-\textfrac{1}{2}\,V\,\xi^1
-\xi^2 
-\left(U -\textfrac{7}{8}\,\kappa-\textfrac{1}{2}\right)\,\xi^3 
\right)
\wedge \theta_2
+\xi^1 \wedge \theta_{12} 
+\xi^2 \wedge \theta_{22}
+\xi^3 \wedge \theta_{23},
\nonumber
\\ \fl
d\theta_3
&=&
\left(\eta_1 -\theta_{22}-\textfrac{1}{2}\,V\,\xi^1-\xi^2
-\left(U-\textfrac{1}{4}\,\kappa\right)\,\xi^3 
\right)
\wedge \theta_3
+\textfrac{1}{2}\,(\kappa+2)\,\xi^2 \wedge \theta_2 
+\xi^1 \wedge \theta_{13}
\nonumber
\\ \fl
&+&
\textfrac{1}{8}\,(\kappa-4)\,
\left(\theta_{22}+\xi^2+\textfrac{1}{8}\,(8\,U-\kappa-4)\,\xi^3\right)\wedge\theta_0
+\xi^2 \wedge \theta_{23}
+\xi^3 \wedge \theta_{12},
\nonumber
\\ \fl
d\xi^1
&=&
-\textfrac{1}{2}\,\left(
\eta_1-3\,\theta_{22}+3\,\xi^2-\left(3\,U-\textfrac{1}{8}\,(\kappa-4)\right)\,\xi^3 
\right)
\wedge \xi^1,
\nonumber
\\ \fl
d\xi^2
&=&
\textfrac{1}{2}\,\left(
\eta_1+\theta_{22}+\textfrac{1}{2}\,V\,\xi^1
+\left(U+\textfrac{1}{8}\,\kappa+\textfrac{3}{2}\right)\,\xi^3 
\right)
\wedge \xi^2
+\left(\textfrac{1}{8}\,(\kappa-4)\,\theta_0 -\theta_3 \right) \wedge \xi^1
-\theta_2 \wedge \xi^3,
\nonumber
\\ \fl
d\xi^3
&=&
\left(
\theta_{22} +\textfrac{1}{2}\,V\,\xi^1 +\xi^2 
\right) 
\wedge \xi^3
-(\kappa+2)\,\left(\theta_2 +\xi^2 \right) \wedge \xi^1,
\nonumber
\\ \fl
d\theta_{11}
&=&
\left(
2\,\eta_1
-3\,\theta_{22}-3\,\xi^2
-\left(3\,U-\textfrac{1}{8}\,(\kappa-4)\right)\,\xi^3 
\right)\wedge \theta_{11}
-\xi^1 \wedge \eta_7
-\xi^2 \wedge \eta_5 
-\xi^3 \wedge \eta_6
\nonumber
\\ \fl
&-&
\left(
\kappa\, \eta_4 
-\textfrac{3}{16}\,V \,\left(3\,\kappa^2+11\,\kappa+4\right)\,\theta_2 
-\kappa\,(U+1)\,\theta_3
-\textfrac{1}{4}\,\kappa\,(\kappa+8)\,\theta_{12} 
-\textfrac{3}{4}\,\kappa\,V\, \theta_{23}
\right.
\nonumber
\\ \fl
&+&
\left.\textfrac{1}{64}\,(\kappa+6)\,(\kappa-4)\,(7\,\kappa+4)\,V\,\xi^2
-\textfrac{5}{64}\,(7\,\kappa+4)\, V^2\,\xi^3
\right)
\wedge \theta_0
+V\,\xi^2 \wedge \theta_{12}
\nonumber
\\ \fl
&+&
\left(
\textfrac{1}{4}\,(5\,\kappa^2+31\,\kappa+24)\,\theta_2 
+(4\,\kappa+3)\,\theta_{23}
+\left((4\,\kappa+3)\,U- \textfrac{1}{4}\,\kappa^2+\textfrac{9}{2}\,\kappa+5\right)\xi^2
\right.
\nonumber
\\ \fl
&+&
\left.
\textfrac{1}{8} \,(15\,\kappa+16)\,V\,\xi^3
\right)
\wedge \theta_1
+\left(\kappa\,\theta_{12} 
+\textfrac{1}{2}\,(\kappa+2)\,V\,\xi^2 
\right)
\wedge \theta_3
\nonumber
\\ \fl
&+&
\left(
(2\,\kappa+3)\,\theta_2
+2\,(\kappa+2)\,\xi^2 
+\textfrac{5}{4}\,V\,\xi^3 
\right)
\wedge \theta_{13},
\label{SE_general}
\\ \fl
d\theta_{12}
&=&
\left(
\eta_1-2\,\theta_{22} -2\,\xi^2-\left(2\,U-\textfrac{1}{2}\,\kappa\right)\,\xi^3 
\right)
\wedge \theta_{12}
-\xi^1 \wedge \eta_5
-\xi^2 \wedge \eta_3
-\xi^3 \wedge \eta_4
\nonumber
\\ \fl
&-&
\left(
\theta_{23} 
+\left((\kappa+1)\,U + \textfrac{3}{8}\,\kappa^2+\textfrac{9}{4}\,\kappa+2\right)\xi^2 
+\textfrac{1}{16}(11\,\kappa+12)\,V\,\xi^3 
\right) 
\wedge \theta_2
\nonumber
\\ \fl
&-&\left(
\textfrac{1}{8}(\kappa-4)\,\theta_0 +\theta_3 
+\textfrac{1}{2}\,V\,\xi^2
\right)\wedge \theta_{22} 
+\left((\kappa+2)\, \xi^2 +\textfrac{3}{4}\,V\,\xi^3\right) \wedge \theta_{23},
\nonumber
\\ \fl
d\theta_{13}
&=&
\textfrac{1}{2}\,\left(
3\,\eta_1 
-5\,\theta_{22} 
-5\,\xi^2 
-\left(5\,U-\textfrac{5}{8}\,\kappa+\textfrac{1}{2}\right)\,\xi^3 
\right)
\wedge \theta_{13}
-\xi^1 \wedge \eta_6
-\xi^2 \wedge \eta_4
-\xi^3 \wedge \eta_5
\nonumber
\\ \fl
&-&
\textfrac{1}{8}\,(\kappa-4)\, 
\left(
\eta_3 
-(\kappa+2)\,\left(U-\textfrac{1}{8}\,\kappa\right)\,\theta_2
-\theta_3 
-\textfrac{1}{8}\,V\,\theta_{22} 
-\theta_{23} 
\right.
\nonumber
\\ \fl
&-&
\left.
   \textfrac{1}{8}\,
    \left((\kappa+2)\,(8\,U-\kappa-4)+6\,V\right)\,
    \xi^2
+\textfrac{1}{64}\,(\kappa+6)\,(7\kappa+4)\,V \,\xi^3
\right) \wedge \theta_0
\nonumber
\\ \fl
&+&
\textfrac{1}{8}\,(\kappa-4)\,\left(
\theta_{22} +\xi^2 +\left(U-\textfrac{1}{8}\,(\kappa+4)\right)\xi^3
\right)\wedge \theta_1
\nonumber
\\ \fl
&-&
\left(
\left(\textfrac{7}{8}\,\kappa^2+4\,\kappa+3\right)\,\theta_3 
+(\kappa+2)\, \theta_{12}
+\textfrac{1}{4}\,(\kappa+2)\,V\,\xi^2 
\right)\wedge \theta_2
\nonumber
\\ \fl
&+&
\left( (2\,\kappa+1)\,\theta_{23} 
+\left(2\,(\kappa+1)\,U+\textfrac{1}{8}\,\kappa^2+\textfrac{11}{4}\,\kappa+3\right)\xi^2
+\textfrac{1}{16}\,(17\,\kappa+20)\,V\,\xi^3
\right)\wedge \theta_3
\nonumber
\\ \fl
&+&\left(\textfrac{3}{2}\,(\kappa+2)\,\xi^2 +V\,\xi^3 \right) \wedge \theta_{12}
+\textfrac{3}{4}\,V\,\xi^2 \wedge \theta_{23},
\nonumber
\\ \fl
d\theta_{22}
&=&
\eta_3 \wedge \xi^1
-\textfrac{1}{2}\,\left(\eta_1+\theta_{22}\right)\wedge \xi^2
+\eta_2 \wedge \xi^3
\nonumber
\\ \fl
d\theta_{23}
&=&
\textfrac{1}{2}\,\left(\eta_1 -3\,\theta_{22}
-3\,\xi^2
-\left(3\,U-\textfrac{11}{8}\,\kappa-\textfrac{1}{4}\right)\xi^3 
\right)
\wedge \theta_{23}
-\xi^1 \wedge \eta_4
-\xi^2 \wedge \eta_2
-\xi^3 \wedge \eta_3
\nonumber
\\ \fl
&-&
\frac{1}{8}\,(\kappa-4)\,\xi^3 \wedge \theta_0 
-\textfrac{1}{8}\,\left(
(3\,\kappa+4)\,\xi^2 
-\left((5\,\kappa-4)\,U+\textfrac{15}{8}\,\kappa^2+12\,\kappa+10\right)
\xi^3 
\right)
\wedge \theta_2
\nonumber
\\ \fl
&+&
\xi^3 \wedge \theta_3
+\textfrac{1}{8}\,
\left(3\,(\kappa+4)\,\theta_2 
+(8\,U-3\,\kappa+4)\,\xi^2
+4\,V\,\xi^3 \right) \wedge \theta_{22},
\nonumber
\\ \fl
d\eta_1
&=&
\textfrac{1}{4}\,\left(\kappa^2+7\,\kappa+4\right) \,\xi^1\wedge \theta_2
+\textfrac{1}{8}\,(\kappa-4)\,\xi^3 \wedge \theta_{22} 
+\kappa\,\xi^1 \wedge \theta_{23} 
+\textfrac{1}{16}\,(7\,\kappa+4)\,V\,\xi^1 \wedge \xi^3
\nonumber
\\ \fl
&+&
\left(\kappa\,U-\textfrac{1}{8}\,\kappa^2+\textfrac{5}{4}\,\kappa+4\right)\,\xi^1 \wedge \xi^2
-\textfrac{1}{8}\,(\kappa-4)\,\xi^2 \wedge \xi^3,
\nonumber
\\ \fl
d\eta_2
&=&
-\textfrac{1}{8}\,\left(
4\,\theta_2 +\left(2\,U+\textfrac{1}{4}\,\kappa+ 3 \right))\,\xi^2
-3\,V\,\xi^3
\right)
\wedge \eta_1
-\left(\theta_{22}+V\,\xi^1+\textfrac{1}{2}\,\xi^2
\right.
\nonumber
\\ \fl
&-&
\left.
(3\,U+2\,\kappa+1)\,\xi^3
\right) \wedge \eta_2
-\textfrac{1}{32} \,\left(\left(32\,U-7\,\kappa^2+38\,\kappa+8\right)\,\xi^1
+96\xi^3\right)
\wedge \eta_3
\nonumber
\\ \fl
&+&\textfrac{1}{16}\,(\kappa-4)\,\theta_0 \wedge 
\left(3\,U + \textfrac{1}{16}\,\left(7\,\kappa^2+40\,\kappa+16\right)\,\xi^1
-3\,\xi^3
\right)
\nonumber
\\ \fl
&-&
\textfrac{1}{2}\,\theta_2 \wedge 
\left(
\theta_{22}
+\left(
\textfrac{1}{128}\,(\kappa+2)
\left(8\,(7\kappa^2+126\kappa+120)\,U+21\kappa^3+198\kappa^2+688\kappa+160\right)
\right.\right.
\nonumber
\\ \fl
&+&
\left.\left.
\textfrac{3}{2}\,(\kappa+3)\,V
\right)\, \xi^1
+\textfrac{3}{8}\,\left(8\,\kappa\,U+3\,\kappa^2+26\,\kappa+24\right)\,\xi^3
\right)
-\textfrac{1}{32}\,\theta_3 \wedge 
\left(
(48\,U+7\,\kappa^2+40\,\kappa
\right.
\nonumber
\\ \fl
&+&
\left.
16)\,\xi^1-48\,\xi^3
\right)
-\theta_{12} \wedge \xi^1
+\textfrac{1}{64}\,\theta_{22} \wedge 
\left(
(7\,\kappa^2+38\,\kappa-8)\,V\,\xi^1
-(16\,U+6\,\kappa+8)\,\xi^2
\right.
\nonumber
\\ \fl
&-&
\left.
(64\,(2\,\kappa+1)\,U+24\,V+128\,(\kappa+1))\,\xi^3
\right)
-\textfrac{1}{32}\,\theta_{23} \wedge 
\left(
(\kappa+2)\,(64\,U+7\,\kappa^2+38\,\kappa
\right.
\nonumber
\\ \fl
&-&
\left.
8)\,\xi^1
-16\,(3\,\kappa+4)\,\xi^3
\right)
+\left(
128\,(\kappa+2)\,U^2+24\,U\,V+4\,(\kappa+2)\,(7\,\kappa^2+70\,\kappa+40)\,U
\right.
\nonumber
\\ \fl
&+&
\left.
(7\,\kappa^2+91\,\kappa+148)\,V+2\,(\kappa+2)\,(7\,\kappa^2+70\,\kappa+24)
\right)\,\xi^1 \wedge \xi^2
+\textfrac{1}{2048}\,V\,\left(
32\,(7\,\kappa^2
\right.
\nonumber
\\ \fl
&+&
\left.
166\,\kappa+152)\,U+2304\,V-49\,\kappa^4-756\,\kappa^3-2324\,\kappa^2
+5152\,\kappa+2752
\right)\,\xi^1 \wedge \xi^3
\nonumber
\\ \fl
&+&
\textfrac{1}{8}\,\left(
3\,V-6\,U^2+2\,(19\,\kappa+16)\,U
+\textfrac{1}{32}\,\left(31\,\kappa^2+1144\,\kappa+1200\right)
\right)\,\xi^2 \wedge \xi^3,
\nonumber
\\ \fl
d \eta_3
&=&
-\textfrac{1}{16}\,\left(
8\,\eta_3-(\kappa-4)\,\theta_0 
+8\,\theta_3+2\,V\,\xi^2
\right) \wedge \eta_1
+\left((\kappa+2)\,(\theta_2 +\xi^2)-\textfrac{1}{2}\,V\,\xi^3\right) \wedge \eta_2
\nonumber
\\ \fl
&-&
\textfrac{1}{32}\,
\left(
48\,\theta_{22}
-24\,V\,\xi^1 
+32\,\xi^2
+(80\,U+ \kappa\,(7\,\kappa+ 36))\,\xi^3 
\right)
\wedge \eta_3 
+(\kappa+2)\,\eta_4 \wedge \xi^1
\nonumber
\\ \fl
&+&
\textfrac{1}{16}\,(\kappa-4)\,\theta_0 \wedge \left(
\theta_{22}
+\left(2\,(\kappa+2)\,(2\,U+1)-V\right)\,\xi^1
+\textfrac{1}{16}\,(48\,U+7\,\kappa^2+40\,\kappa+16)\xi^3 
\right)
\nonumber
\\ \fl
&+&
\textfrac{1}{8}\,(\kappa-4)\,\theta_1 \wedge \xi^1
+\textfrac{1}{32}\,\theta_2 \wedge \left(
V\,(80\,(\kappa+2)\,U+7\,\kappa^3+36\,\kappa^2+10\,\kappa-56)\,\xi^1
\right.
\nonumber
\\ \fl
&+&
%\left.
4\,(\kappa^2+7\,\kappa+4)\,\xi^2
-\left(
(\kappa+2)\,(7\,\kappa^2+126\,\kappa+120)\,U
+24\,(\kappa+3)\,V
\right.
\nonumber
\\ \fl
&+&
\left.
\left.
\textfrac{1}{8}\,(\kappa+2)\,(21\,\kappa^3+198\,\kappa^2+688\,\kappa+160)
\right)\,\xi^3
\right)
-\textfrac{1}{2}\,\theta_3 \wedge \left( 
\theta_{22}
+2\,((\kappa+2)\,(2\,U+1)
\right.
\nonumber
\\ \fl
&-&
\left.
V)\,\xi^1
+\left(6\,U+\textfrac{7}{16}\,\kappa^2+\textfrac{5}{2}\,\kappa+1\right)\xi^3
\right)
+
\left(
(\kappa+2)\,\left(U+\textfrac{3}{8}\,(\kappa+4)\right)\,\xi^1
+\xi^3
\right) \wedge \theta_{12}
\nonumber
\\ \fl
&-&\theta_{13} \wedge \xi^1
+\textfrac{1}{64}\,V\,\theta_{22} \wedge \left(
8\,\xi^2-\left(7\,\kappa^2+38\,\kappa-8\right)\,\xi^3
\right)
-\textfrac{1}{4}\,\theta_{23} \wedge \left(
\kappa\,(V\,\xi^1-2\,\xi^2)
\right.
\nonumber
\\ \fl
&+&
\left.
\textfrac{1}{8}\,(\kappa+2)\,(64\,U +7\,\kappa^2+38\,\kappa-8)
\,\xi^3
\right)
-\textfrac{1}{64}\,V\,
\left(32\,(11\,\kappa+5)\,U+8\,V+21\,\kappa^3
\right.
\nonumber
\\ \fl
&+&
\left.
124\,\kappa^2+44\,\kappa-208
\right)\,\xi^1 \wedge \xi^2
-\textfrac{1}{256}\,V^2\,\left(160\,U+21\,\kappa^2+210\,\kappa+104\right)\xi^1 \wedge \xi^3
\nonumber
\\ \fl
&-&
\textfrac{1}{64}\,\left(
128\,(\kappa+2)\,U^2+24\,U\,V
+4\,(\kappa+2)\,(7\,\kappa^2+70\,\kappa+40)\,U
+(7\,\kappa^2 + 105\,\kappa+ 156)\,V
\right.
\nonumber
\\ \fl
&+&
\left.
2\,(\kappa+2)\,(7\,\kappa^2+70\,\kappa+24)
\right)\xi^2 \wedge \xi^3,
\nonumber
\\ \fl
d\eta_4
&=&
\eta_8 \wedge \xi^1
+\left(\theta_3 - \textfrac{1}{8}\,(\kappa-4)\,\theta_0
-\textfrac{3}{4}\,V\,\xi^2\right) \wedge \eta_2
+\textfrac{1}{32}\,\left(
8\,(11\,\kappa+7)\,\theta_2
+48\,\theta_{23}
\right.
\nonumber
\\ \fl
&-&
\left.
(7\,\kappa^2+18\,\kappa-88)\,\xi^2
+56\,\xi^3
\right) \wedge \eta_3
+\left(\eta_1 -3\,\theta_{22} -3\,\xi^2 
-\left(3\,U-\textfrac{1}{4}\,\kappa+2\right)\,\xi^3
\right) \wedge \eta_4
\nonumber
\\ \fl
&+&
\textfrac{1}{256}\,(\kappa-4)\,\left(
4\,(11\,\kappa+20)\,\theta_2
-4\,(8\,U-3\,\kappa+4)\,\theta_{22}
+48\,\theta_{23}
-(48\,U+7\,\kappa^2+8\,\kappa-48)\,\xi^2
\right.
\nonumber
\\ \fl
&-&
\left.
32\,((\kappa+2)\,(2\,U+1)-V)\,\xi^3
\right)
 \wedge \theta_0
+\textfrac{1}{512}\,
\left(
64\,(11\,\kappa+20)\,\theta_3 
+16\,V\,(13\,\kappa+20)\,\theta_{22}
\right.
\nonumber
\\ \fl
&-&
\left. 
32\,(24\,(\kappa+2)\,U-13\,\kappa^2-40\,\kappa-16)\,\theta_{23}
+2\,\left(
8\,(\kappa+2)\,(7\,\kappa^2+106\,\kappa+136)\,U
\right.\right.
\nonumber
\\ \fl
&+&
\left.\left.
192\,(\kappa+3)\,V
+(\kappa+2)\,(21\,\kappa^3+138\,\kappa^2+304\,\kappa-160)
\right)\,\xi^2
-V\,\left(
16\,(65\,\kappa+172)\,U
\right.\right.
\nonumber
\\ \fl
&+&
\left.\left.
77\,\kappa^3+444\,\kappa^2+544\,\kappa-768
\right)\,\xi^3
\right) \wedge \theta_2
+
\left(
\left(U- \textfrac{1}{8}\,(3\,\kappa-4)\right)\,\theta_{22} 
-\textfrac{3}{2}\,\theta_{23}
\right.
\nonumber
\\ \fl
&+&
\left.
\textfrac{1}{32}\,\left(48\,U +7\,\kappa^2+8\,\kappa-48\right)\,\xi^2 
+((\kappa+2)\,(2\,U+1)-V)\,\xi^3
\right) \wedge \theta_3
\nonumber
\\ \fl
&+&\textfrac{1}{64}\,\left(
24\,(\kappa+4)\,(\theta_{22}+\xi^2) 
+(8\,(3\,\kappa+20)\,U+9\,\kappa^2+48\,\kappa+112)\,\xi^3
\right)
\wedge \theta_{12}
\nonumber
\\ \fl
&+&
\textfrac{1}{32}\,V\,
\left(
(24\,U+7\,\kappa^2+25\,\kappa-44)\,\xi^2
-8\,V\,\xi^3
\right) \wedge \theta_{22}
+\textfrac{1}{128}\,\left(
4\,(32\,(3\,\kappa+7)\,U+12\,V
\right.
\nonumber
\\ \fl
&+&
\left.
7\,\kappa^3+32\,\kappa^2-4\,\kappa-48)\,\xi^2
-V\,(96\,U+21\,\kappa^2+154\,\kappa+8)\,\xi^3
\right)\wedge \theta_{23}
\nonumber
\\ \fl
&+&
\textfrac{1}{2048}\,V\,\left(
32\,(7\,\kappa^2+236\,\kappa+504)\,U+2560\,V
-49\,\kappa^4-84\,\kappa^3+1644\,\kappa^2+6560\,\kappa
\right.
\nonumber
\\ \fl
&-&
\left.
3904
\right)\,\xi^2 \wedge \xi^3,
\nonumber
\\ \fl
d\eta_5
&=&
\eta_9 \wedge \xi^1
+\eta_8 \wedge \xi^3
-\left(
\textfrac{1}{4}\,(\kappa-4)\,\theta_0
-2\,\theta_3
-2\,\theta_{12}
-\textfrac{3}{2}\,V\,\xi^2
\right) \wedge \eta_3
+
\left((\kappa+3)\,\theta_2
\right.
\nonumber
\\ \fl
&+&
\left.
(\kappa+2)\,\xi^2 
+\textfrac{5}{4}\,V\,\xi^3
\right) \wedge \eta_4
+\textfrac{1}{16}\,\left(
24\,\eta_1 - 56\,(\theta_{22}+\xi^2)
-(56\,U -9\,\kappa+8)\,\xi^3
\right)\wedge \eta_5
\nonumber
\\ \fl
&+&
\textfrac{1}{64}\,(\kappa-4)\,\theta_0 \wedge \left(
(8\,(\kappa+1)\,U +3\,\kappa^2+18\,\kappa+16)\,\theta_2
-16\,\theta_{12}
+4\,V\,\theta_{22}
+8\,(\kappa+2)\,\theta_{23}
\right.
\nonumber
\\ \fl
&+&
\left.		
4\,(2\,(\kappa+2)\,(2\,U+1)-V)\,\xi^2
\right)
+\textfrac{1}{8}\,(\kappa-4)\,\theta_1 \wedge (\theta_{22}+\xi^2)
\nonumber
\\ \fl
&+&\textfrac{1}{64}\,\theta_2 \wedge \left(
8\,\left(8\,(\kappa+1)\,U
+
3\,\kappa^2+18\,\kappa-16\right)\,\theta_3
+16\,\left(8\,(\kappa+2)-3\,\kappa^2-11\,\kappa-4\right)\,\theta_{12}
\right.
\nonumber
\\ \fl
&-&
\left.
16\,V\,(3\,\kappa+5)\,\theta_{23}
+V\,\left(32\,(3\,\kappa+8)\,U+7\,\kappa^3+71\,\kappa^2+126\,\kappa-8\right)\,\xi^2
\right.
\nonumber
\\ \fl
&+&
\left.
V^2\,(11\,\kappa+12)\,\xi^3 
\right)
+\theta_3 \wedge \left(
+2\,\theta_{12}-(\kappa+2)\theta_{23} 
-\left((\kappa+2)\,(2\,U+1)-V\right)\,\xi^2
\right)
\nonumber
\\ \fl
&+&
\theta_{12} \wedge \left(
(\kappa-1)\,\theta_{23}
-\left((\kappa+5)\,U+\textfrac{1}{2}\,V-\textfrac{3}{8}\,\kappa\,(\kappa+6)\right) \xi^2
-\textfrac{1}{32}\,V\,\left(32\,U +7\,\kappa^2+82\,\kappa
\right.\right.
\nonumber
\\ \fl
&+&
\left.\left.
24\right)\,
\xi^3
\right)
-\theta_{13} \wedge (\theta_{22}+\xi^2)
-\textfrac{1}{8}\,V^2\,\theta_{22}\wedge\xi^2
+\textfrac{1}{8}\,V\,\theta_{23} \wedge \left(
2,(3\,\kappa+8)\, \xi^2 
+3\,V\,\xi^3
\right)
\nonumber
\\ \fl
&+&\textfrac{1}{256}\,V^2\left(
160\,U+21\,\kappa^2+210\,\kappa+104
\right)\,\xi^2 \wedge \xi^3,
\nonumber
\\ \fl
d\eta_6
&=&
\eta_{10} \wedge \xi^1
+\eta_8 \wedge \xi^2
+\eta_9 \wedge \xi^3
+\textfrac{1}{32}\,\eta_3 \wedge \left(
(\kappa-4)\,(5\,V\,\theta_0+8\,\theta_1)
-40\,\theta_{13}
\right)
\nonumber
\\ \fl
&+&
\eta_4 \wedge \left(
2\,\kappa\,\theta_3
-\textfrac{1}{8}\,\kappa\,(\kappa-4)\,\theta_0
+ V \, \xi^2
\right)
-\textfrac{1}{2}\,\eta_5 \wedge \left(
(\kappa+2)\,(4\,\theta_2 + 5\,\xi^2) + 3\,V\,\xi^3
\right)
\nonumber
\\ \fl
&-&
\eta_6 \wedge \left(\
2\,(\eta_1+2\,\theta_{22}+2\,\xi^2) 
+\left(4\,U-\textfrac{1}{8}\,(3\,\kappa-4)\right)\,\xi^3
\right)
+\textfrac{1}{2048}\,(\kappa-4)\,\theta_0 \wedge \left(
64\,(\kappa-4)\,\theta_1
\right.
\nonumber
\\ \fl
&+&
\left.
8\,(40\,(\kappa+2)\,U-19\,\kappa^2-72\,\kappa-56)\,\theta_2
+128\,(\kappa+2)\,\theta_{12}
+640\,\theta_{13}
-64\,V\,\theta_{23}
\right.
\nonumber
\\ \fl
&+&
\left.
4\,V\,(16\,(7\,\kappa+16)\,U+20\,V+7\,\kappa^3+10\,\kappa^2-152\,\kappa-256)\,\xi^2
-V^2\,(160\,U+35\,\kappa^2
\right.
\nonumber
\\ \fl
&+&
\left.
302\,\kappa+152)\,\xi^3
\right)
+\textfrac{1}{32}\,\theta_1 \wedge \left(
4\,(256\,\kappa^3-3\,\kappa-4-2\,(\kappa-4)\,(\kappa+2)\,U)\,\theta_2 
\right.
\nonumber
\\ \fl
&+&
\left.
(\kappa-4)\,(8\,\theta_3-V\,\theta_{22}
+4\,\theta_{23}
+((\kappa+2)\,(8\,U-\kappa-4)+V)\,\xi^2
+(40\,U+7\,\kappa^2
\right.
\nonumber
\\ \fl
&+&
\left.
41\,\kappa+4)\,\xi^3)
\right)
-\theta_2 \wedge \left(
\textfrac{1}{32}\,V\,(13\,\kappa^2+4\,\kappa-8)\,\theta_3
+\textfrac{1}{16}\,(40\,(\kappa+2)\,U-25\,\kappa^2-112\,\kappa
\right.
\nonumber
\\ \fl
&-&
\left.
64)\,\theta_{13}
\right)
-\textfrac{1}{32}\,\theta_3 \wedge \left(
((\kappa-4)\,(4\,\kappa\,U+5\,V)+2\,\kappa\,(4\,\kappa-1)
 )\,\theta_0 
-4\,\kappa\,(7\,\kappa +20)\,\theta_{12}
\right.
\nonumber
\\ \fl
&+&
\left.
8\,V\,(4\,\kappa +5)\,\theta_{23}
+V\,(80\,(\kappa+1)\,U+7\,\kappa^3+12\,\kappa^2-56\,\kappa-64)\,\xi^2
\right)
\nonumber
\\ \fl
&+&
\textfrac{1}{64}\,(\kappa-4)\,\theta_{11} \wedge \left(
8\,(\theta_{22}+\xi^2)
+(8\,U-\kappa-4)\,\xi^3
\right)
+ \textfrac{1}{8}V\,\theta_{12} \wedge \left(
13\,(\kappa+2)\,\xi^2 + 4\,V\,\xi^3 
\right)
\nonumber
\\ \fl
&+&
\theta_{13} \wedge \left(
\textfrac{1}{2}\,(7\,\kappa+2)\,\theta_{23}
+\textfrac{1}{8}\,(8\,(\kappa-3)\,U - 5\,V + 8\,(5\,\kappa+4))\,\xi^2
+\textfrac{1}{128}\,V\,(160\,U+35\,\kappa^2
\right.
\nonumber
\\ \fl
&+&
\left.
430\,\kappa+152)\,\xi^3
\right)
+\textfrac{3}{16}\,V^2\,\theta_{23} \wedge \xi^2,
\nonumber
\\ \fl
d\eta_7
&=&
\eta_{11} \wedge \xi^1
+\eta_9 \wedge \xi^2
+\eta_{10} \wedge \xi^3
-\kappa\,\eta_8 \wedge \theta_0
-3\,\eta_3 \wedge \theta_{11} 
+\eta_4 \wedge \left(
\textfrac{3}{4}\,\kappa\,V\,\theta_0+(5\,\kappa+3)\,\theta_1
\right)
\nonumber
\\ \fl
&+&
\textfrac{1}{8}\,\eta_5 \wedge \left(
(2\,\kappa^2+17\,\kappa-4)\,\theta_0+8\,(\kappa-1)\,\theta_3-10\,V\,\xi^2
\right)
-\eta_6 \wedge \left(
(3\,\kappa+5)\,\theta_2+3\,(\kappa+2)\, \xi^2
\right.
\nonumber
\\ \fl
&+&
\left.
\textfrac{7}{4}\,V\,\xi^3
\right)
-\textfrac{1}{16}\,\eta_7 \wedge \left(
40\,\eta_1-72\,(\theta_{22}+\xi^2)-(72\,U-3\,(\kappa-4))\,\xi^3
\right)
\nonumber
\\ \fl
&+&
\textfrac{1}{256}\,\theta_0 \wedge \left(
8\,(\kappa-4)\,(4\,(4\,\kappa+3)\,U-\kappa^2+18\,\kappa-20)\,\theta_1
+4\,V^2\,(26\,\kappa^2+57\,\kappa+28)\,\theta_2
\right.
\nonumber
\\ \fl
&-&
\left.
8\,V\,(48\,\kappa\,U+7\,\kappa^3+26\,\kappa^2-64\,\kappa-32)\,\theta_3
-96\,(\kappa-4)\,\theta_{11}
-16\,V\,(9\,\kappa^2+35\,\kappa+4)\,\theta_{12}
\right.
\nonumber
\\ \fl
&+&
\left.
64\,(4\,\kappa\,U+\kappa^2+2\,\kappa-8)\,\theta_{13}
-(7\,\kappa+4)\,(\kappa^2-18\,\kappa-64)\,V^2\,\xi^2
-10\,(7\,\kappa+4)\,V^3\, \xi^3
\right)
\nonumber
\\ \fl
&+&
\textfrac{1}{64}\,\theta_1 \wedge \left(
4\,(\kappa-1)\,(19\,\kappa+20)\,V\,\theta_2
+32\,(2\,(5\,\kappa+3)\,U-(8\,\kappa+5)\,(\kappa-2))\,\theta_3
\right.
\nonumber
\\ \fl
&+&
\left.
16\,(6\,\kappa^2+39\,\kappa+24)\,\theta_{12}
-144\,(\kappa+1)\,V\,\theta_{23}
-V\,(96\,(4\,\kappa+3)\,U+35\,\kappa^3+63\,\kappa^2
\right.
\nonumber
\\ \fl
&-&
\left.
302\,\kappa-216)\,\xi^2
+(5\,\kappa-12)\,V^2\,\xi^3
\right)
+\textfrac{1}{8}\,\theta_2 \wedge \left(
(24\,(\kappa+2)\,U-15\,\kappa^2-88\,\kappa-56)\,\theta_{11}
\right.
\nonumber
\\ \fl
&-&
\left.
2\,(3\,\kappa+7)\,V\, \theta_{13}
\right)
+\textfrac{1}{8}\,\theta_3 \wedge \left(
24\,\theta_{11}-(\kappa+2)\,(4\,V\,\theta_{12}+16\,\theta_{13}+V^2\,\xi^2)
\right)
\nonumber
\\ \fl
&+&
\theta_{11} \wedge \left(
3\,(\kappa+2)\,\theta_{23}
+\textfrac{1}{8}\,(24\,(\kappa-1)\,U-6\,V-3\,\kappa^2+54\,\kappa+48)\,\xi^2
+\textfrac{3}{64}\,(32\,U+7\,\kappa^2
\right.
\nonumber
\\ \fl
&+&
\left.
98\,\kappa+40)\,\xi^3
\right)
-\theta_{12} \wedge \left(
3\,(\kappa+1)\,\theta_{13}
-\textfrac{1}{4}\,V^2\,\xi^2
\right)
+\textfrac{1}{8}\,V\,\theta_{13} \wedge \left(
18\,(\kappa+2)\,\xi^2
+5\,V\,\xi^3 
\right).
\nonumber
\end{eqnarray}
}

%==========================================================================

The structure equations of the symmetry pseudo-group of Eq. (\ref{main}) in the case of  
$\kappa  = -1$ have the form
%----------------------------------------------------------------------------------
%       Exceptional case
%----------------------------------------------------------------------------------
{\small
\begin{eqnarray}
\fl  
d\theta_0
&=&
\eta_1 \wedge \theta_0+\xi^1 \wedge \theta_1+\xi^2 \wedge \theta_2
+\xi^3 \wedge \theta_3,
\nonumber
\\ \fl
d\theta_1
&=&
\textfrac{3}{2}\,\eta_1 \wedge \theta_1
+2\,\eta_2 \wedge \theta_3
+\eta_3 \wedge \theta_0
-\textfrac{1}{8}\,\left(
12\,\theta_{22}+16\,U\,\xi^1+12\,\xi^2+7\,\xi^3
\right) \wedge \theta_1
+\xi^1 \wedge \theta_{11}
\nonumber
\\ \fl
&+&
\xi^2 \wedge \theta_{12}
+\xi^3 \wedge \theta_{13},
\nonumber
\\ \fl
d\theta_2
 &=&
\textfrac{1}{8}\,\left(
4\,(\eta_1-\theta_{22} -\xi^2)-3\,\xi^3
\right)\wedge \theta_2
+\xi^1 \wedge \theta_{12}+\xi^2 \wedge \theta_{22}+\xi^3 \wedge \theta_{23},
\nonumber
\\ \fl
d\theta_3
&=&
\left(
\eta_1-\theta_{22}-U\,\xi^1-\xi^2-\textfrac{5}{8}\,\xi^3 
\right)\wedge \theta_3
+\eta_2 \wedge \theta_2
-\textfrac{5}{8}\,(\theta_{22}+\xi^2) \wedge \theta_0 
+\xi^1 \wedge \theta_{13}
\nonumber
\\ \fl
&+&
\xi^2 \wedge \theta_{23}
+\xi^3 \wedge \theta_{12},
\nonumber
\\ \fl
d\xi^1
&=&
\textfrac{1}{8}\,\left(
12\,\theta_{22}-4\,\eta_1+12\,\xi^2+7\,\xi^3
\right) \wedge \xi^1,
\nonumber
\\ \fl
d\xi^2
&=&
-\textfrac{1}{8}\,\left(
5\,\theta_0 + 8\,\theta_3\right) \wedge \xi^1
+\textfrac{1}{8}\,\left(
4\,(\eta_1 +\theta_{22})+3\,\xi^3
\right) \wedge \xi^2
-\left(\eta_2 +\theta_2 \right)\wedge \xi^3,
\nonumber
\\ \fl
d\xi^3
&=&
-(2\,\eta_2 +\theta_2) \wedge \xi^1
+(\theta_{22} 
+U\,\xi^1+\xi^2) \wedge \xi^3,
\nonumber
\\ \fl
d\theta_{11}
&=&
2\,\eta_1 \wedge \theta_{11}
+\eta_2 \wedge \left(\textfrac{1}{2}\,U\,\theta_0+4\,\theta_{13}\right)
+6\,\eta_3 \wedge \theta_1
+\eta_4 \wedge \theta_0
+\eta_7 \wedge \xi^1
+\eta_5 \wedge \xi^2
+\eta_6 \wedge \xi^3
\nonumber
\\ \fl
&+&
\textfrac{1}{8}\,\left( 
 30\,\theta_1
-24\,U\,\theta_2 
-11\,\theta_3 
-14\,\theta_{12}
-24\,U\,\theta_{23}
-22\,U\,\xi^2 
+15\,U^2\,\xi^3 
\right)
\wedge \theta_0
\nonumber
\\ \fl
&+&
\textfrac{5}{8}\,\left(
10\,\theta_2 
+8\,\theta_{23} 
+7\,\xi^2 
-4\,U\,\xi^3 
\right)
\wedge \theta_1
+(2\,U\,\theta_3 -\theta_{13})\wedge \theta_2
-(\theta_{12}+U\,\xi^2) \wedge \theta_3
\nonumber
\\ \fl
&-&
\left(3\,(\theta_{22}+\xi^2)+\textfrac{7}{4}\,\xi^3\right) \wedge \theta_{11}
+2\,U\,\xi^2 \wedge \theta_{12}
+3\,U\,\xi^3 \wedge \theta_{13},
\label{SE_exceptional}
\\ \fl
d\theta_{12}
&=&
\eta_1 \wedge \theta_{12}
+\eta_2 \wedge \left(2\,\theta_{23}+\textfrac{3}{2}\,\xi^2\right)
+\eta_3 \wedge (\theta_2+\xi^2)
+\eta_4 \wedge\xi^3
+\eta_5 \wedge \xi^1
-\textfrac{5}{8}\,\theta_0 \wedge \left(\theta_2-\theta_{22}\right)
\nonumber
\\ \fl
&+&
\textfrac{1}{4}\,\theta_2 \wedge \left(7\,\xi^2 +2\,U\,\xi^3\right)
+\theta_3 \wedge (\theta_{22}+\xi^2)
+\theta_{12} \wedge \left(2\,(\theta_{22}+\xi^2)+\textfrac{5}{4}\,\xi^3\right)
\nonumber
\\ \fl
&+&
\theta_{23}\wedge\left(2\,\xi^2 -U\, \xi^3\right)
+\textfrac{1}{4}\,U\,\xi^2 \wedge \xi^3,
\nonumber
\\ \fl
d\theta_{13}
&=&
\textfrac{3}{2}\,\eta_1 \wedge \theta_{13}
+3\, \eta_2 \wedge \left(\textfrac{5}{16}\,\theta_0+\theta_{12}\right)
+\eta_3 \wedge \left(
\textfrac{5}{8}\,\theta_0
+3\,\theta_3 \right)
+\eta_4 \wedge \xi^2
+\eta_5 \wedge \xi^3
+\eta_6 \wedge \xi^1
\nonumber
\\ \fl
&-&
\textfrac{5}{64}\,\theta_0 \wedge \left(
17\,\theta_2 
+24\theta_3
-8\,U\,\theta_{22} 
+8\,\theta_{23}
-(8\,U - 14)\xi^2 
-2\,U\,\xi^3 
\right)
+\textfrac{5}{8}\,\theta_1 \wedge (\theta_{22}+\xi^2) 
\nonumber
\\ \fl
&+&
\textfrac{1}{4}\,\theta_2 \wedge \left(
13\,\theta_3 
+4\,\theta_{12} 
+2\,U\,\xi^2
\right)
-\textfrac{1}{8}\,\theta_3 \wedge \left(
16\,\theta_{23}+21\,\xi^2+12\,U\,\xi^3
\right)
+2\,U\,\xi^3 \wedge \theta_{12} 
\nonumber
\\ \fl
&+&
\textfrac{1}{2}\,\theta_{13} \wedge \left(5\,(\theta_{22}+\xi^2)+3\,\xi^3\right) 
+U\,\xi^2 \wedge \theta_{23},
\nonumber
\\ \fl
d\theta_{22}
&=&
-\textfrac{1}{2}\,\eta_1 \wedge \xi^2
+\textfrac{3}{2}\,\eta_2 \wedge \xi^1
+\eta_2 \wedge \xi^3 
+\eta_3 \wedge \xi^1 
+\theta_2 \wedge \left(\textfrac{21}{8}\,\xi^1 +\xi^3\right)
+\theta_3 \wedge \xi^1
\nonumber
\\ \fl
&-&
\textfrac{1}{8}\,\theta_{22} \wedge \left( 
8\,\xi^2 +3\,\xi^3
\right)
+2\,\theta_{23} \wedge \xi^1 
-\textfrac{1}{4}\,\xi^1 \wedge 
(7\,\xi^2-U\,\xi^3),
\nonumber
\\ \fl
d\theta_{23}
&=&
\textfrac{1}{2}\,\eta_1 \wedge \theta_{23}
+\eta_2 \wedge \left(\theta_{22}+\xi^2+\textfrac{3}{2}\,\xi^3 \right)
+\eta_3 \wedge \xi^3
+\eta_4 \wedge\xi^1
-\textfrac{5}{8}\,\theta_0 \wedge \xi^3 
\nonumber
\\ \fl
&+&
\textfrac{1}{16}\,\theta_2 \wedge \left(
72\,(\theta_{22}+\xi^2)+49\,\xi^3
\right)
-\textfrac{3}{8}\,\theta_{22}\,\left(4\,\theta_{23}+3\,\xi^2\right)
+\textfrac{3}{2}\,\theta_{23} \wedge \left(\xi^2 +2\,\xi^3 \right)
+\textfrac{77}{32}\,\xi^2 \wedge \xi^3,
\nonumber
\\ \fl
d\eta_1
&=&
-\eta_3 \wedge \xi^1
+\textfrac{5}{8}\,\left((\theta_0-\theta_2) \wedge \xi^1
+(\theta_{22}+\xi^2) \wedge \xi^3 
\right),
\nonumber
\\ \fl
d\eta_2
&=&
-\textfrac{1}{4}\,\eta_2 \wedge 
\left(
2\,\eta_1-2\,\theta_{22}-4\,U\,\xi^1-2\,\xi^2-\xi^3
\right)
-\eta_3 \wedge \xi^3
-\textfrac{5}{16}\,\theta_0 \wedge (\xi_1-2\,\xi^3)
\nonumber
\\ \fl
&+&
\theta_2 \wedge \left(U\,\xi^1-\textfrac{13}{8}\,\xi^3\right)
-\textfrac{1}{2}\,\theta_3 \wedge \xi^1
+\textfrac{1}{2}\,\theta_{22} \wedge \xi^2
-\theta_{23} \wedge \xi^3
-\textfrac{1}{16}\,(8\,U\,\xi^1 -21\,\xi^3) \wedge \xi^2,
\nonumber
\\ \fl
d\eta_3
&=&
\textfrac{1}{16}\,\eta_1 \wedge \left( 
5\,\theta_0+ 8\,\eta_3 
\right)
-\textfrac{1}{16}\,\eta_2 \wedge \left(
20\,\theta_{22}+8\,U\,\xi^1+20\,\xi^2+15\,\xi^3 
\right)
-\eta_4 \wedge \xi^1 
-\textfrac{5}{8}\,\theta_1 \wedge \xi^1
\nonumber
\\ \fl
&+&
\textfrac{1}{4}\,\eta_3 \wedge \left(
6\,\theta_{22}+8\,U\,\xi^1+6\,\xi^2+\xi^3
\right)
-\textfrac{5}{32}\,\theta_0 \wedge \left(
6\,\theta_{22}+(4\,U-3)\,\xi^1+\frac{15}{16}\,\xi^2 \wedge \theta_0+\xi^3
\right)
\nonumber
\\ \fl
&+&\textfrac{1}{64}\,\theta_2 \wedge \left(
192\,U\,\xi^1-40\,\xi^2-85\,\xi^3
\right)
+\textfrac{1}{8}\,\theta_3 \wedge \left(
6\,\xi^1-5\,\xi^3
\right)
+\textfrac{9}{8}\,\theta_{12} \wedge \xi^1
\nonumber
\\ \fl
&+&
\textfrac{5}{8}\,\theta_{22} \wedge (\xi^2+U\,\xi^3)
+\theta_{23} \wedge \left(
3\,U\,\xi^1-\textfrac{5}{8}\,\xi^3
\right)
-\textfrac{1}{8}\,U\,\xi^1 \wedge 
(22\,\xi^2+15\,U\,\xi^3)
\nonumber
\\ \fl
&+&\textfrac{5}{32}\,(4\,U -7)\xi^2 \wedge \xi^3,
\nonumber
\\ \fl
d\eta_4
&=&
\eta_8 \wedge \xi^1
+\textfrac{1}{16}\,\eta_2 \wedge 
\left(
+48\,\eta_3
-30\,\theta_0
+89\,\theta_2
+16\,U\,\theta_{22}
+92\,\theta_{23}
+2\,(8\,U+39)\,\xi^2
-12\,U\,\xi^3
\right)
\nonumber
\\ \fl
&-&
\textfrac{1}{8}\,\eta_3 \wedge \left(
17\,\theta_2 
+16\,\theta_{23}
+9\,\xi^2
-8\,U\,\xi^3
\right)
-\eta_4 \wedge \left(
\eta_1 
-3\,\theta_{22}
-3\,\xi^2-\textfrac{23}{8}\,\xi^3\right)
%%%
\nonumber
\\ \fl
&+&
\textfrac{5}{64}\,\theta_0 \wedge 
\left(
17\,\theta_2 
-13\,\theta_{22}
+16\,\theta_{23}
-4\,\xi^2
-\left(8\,U+\textfrac{77}{4}\right)\xi^3 
\right)
%%%
+\textfrac{1}{64}\,\theta_2 \wedge \left(
64\,(U\,\theta_{22}+\theta_{23})
\right.
\nonumber
\\ \fl
&+&
\left.
(64\,U+91)\xi^2 
+74\,U\,\xi^3
\right)
%%%
-\textfrac{1}{32}\,\theta_3 \wedge \left(
26\,(\theta_{22}+\xi^2)+77\,\xi^3
\right)
-\textfrac{1}{16}\,\theta_{12} \wedge \left(
18\,(\theta_{22}+\xi^2)
\right.
\nonumber
\\ \fl
&+&
\left.
31\xi^3
\right)
-\textfrac{1}{2}\,U\,\theta_{22} \wedge \xi^2
+\textfrac{3}{8}\,\theta_{23} \wedge \left(\xi^2
+U\,\xi^3
\right)
+\textfrac{13}{16}\,U\,\xi^2 \wedge \xi^3,
\nonumber
\\ \fl
d\eta_5
&=&
\eta_9 \wedge \xi^1
+\eta_8 \wedge \xi^3
-\textfrac{1}{8}\,\eta_2 \wedge \left(
15\,\theta_0
-4\,U\,\theta_2 
+24\,\theta_3
+ 4\,\theta_{12}
-12\,U\,\xi^2
\right)
\nonumber
\\ \fl
&-&
\textfrac{1}{4}\,\eta_3 \wedge \left(
5\,\theta_0 
+8\,\theta_3 
-8\,U\,\theta_2
+16\,\theta_{12}
-8\,U\,\xi^2 \right)
-\eta_4 \wedge \left(
4\,\eta_2 +2\,\theta_2 -\xi^2 +2\,U\,\xi^3
\right)
\nonumber
\\ \fl
&-&
\textfrac{1}{8}\,\eta_5 \wedge \left(
12\,\eta_1
-28\,(\theta_{22}+\xi^2)
-17\,\xi^3
\right)
%%%%%%%%%%%
-\textfrac{5}{64}\,\theta_0 \wedge \left(
(16\,U-41)\theta_2 
-16\,\theta_3
-32\,\theta_{12} 
\right.
\nonumber
\\ \fl
&-&
\left.
24\,\theta_{23}
+2\,(8\,U-7)\xi^2 
+4\,U\,\xi^3
\right)
-\textfrac{5}{8}\,\theta_1 \wedge (\theta_{22}+\xi^2)
-\textfrac{1}{8}\,\theta_2 \wedge \left(
41\,\theta_3
+46\,\theta_{12}
\right.
\nonumber
\\ \fl
&-&
\left.
16\,U\,\xi^2
-11\,U^2\,\xi^3
\right)
+\textfrac{1}{4}\,\theta_3 \wedge \left(
4\,U\,\theta_{22}+12\,\theta_{23}
+(4\,U-7)\,\xi^2 \wedge \theta_3
-2\,U\,\xi^3
\right)
\nonumber
\\ \fl
&+&
\textfrac{1}{8}\, \theta_{12} \wedge \left(
32\,\theta_{23}+23\,\xi^2-18\,U\,\xi^3
\right)
-\theta_{13} \wedge (\theta_{22}+\xi^2) 
+U\,\theta_{23} \wedge \left(
2\,(\theta_2+\xi^2)+U\,\xi^3
\right)
\nonumber
\\ \fl
&+&
\textfrac{13}{8}\,U^2\,\xi^2 \wedge \xi^3,
\nonumber
\\ \fl
d\eta_6
&=&
\eta_{10} \wedge \xi^1
+\eta_8 \wedge \xi^2 
+\eta_9 \wedge \xi^3
-\textfrac{1}{8}\,\eta_2 \wedge \left(
15\,\theta_1
-12\,U\,\theta_3
+8\,U\,\theta_{12} 
-6\,\theta_{13}
\right)
\nonumber
\\ \fl
&+&
\textfrac{1}{8}\,\eta_3 \wedge 
\left(
5\,U\,\theta_0 
-10\,\theta_1
+72\,U\,\theta_3
-56\,\theta_{13} 
\right)
-\textfrac{1}{8}\,\eta_4 \wedge \left(
5\,\theta_0
+16\,\theta_3
+8\,U\,\xi^2
\right)
\nonumber
\\ \fl
&-&
\eta_5 \wedge \left(
+5\,\eta_2 +2\,\theta_2+3\,U\,\xi^3
\right)
-\eta_6 \wedge \left(
2\,\eta_1-4\,\theta_{22}-4\,\xi^2 -\textfrac{19}{8}\,\xi^3 \right)
%%%%
\nonumber
\\ \fl
&+&
\textfrac{5}{64}\,\theta_0\wedge\left(
+10\,\theta_1
-12\,U\,\theta_2
-(72\,U-11)\theta_3
-4\,\theta_{12} 
+56\,\theta_{13}
-16\,U\,\theta_{23}
-14\,U\,\xi^2
\right.
\nonumber
\\ \fl
&+&
\left.
11\,U^2\,\xi^3
\right)
+\textfrac{5}{32}\,\theta_1 \wedge \left(
19\,\theta_2
+4\,U\,\theta_{22}
+12\,\theta_{23}
+2\,(4\,U+7)\,\xi^2 
-2\,U\,\xi^3
\right)
\nonumber
\\ \fl
&+&
\textfrac{1}{4}\,\theta_2 \wedge \left(
53\,U\,\theta_3
+4\,U\,\theta_{12}
-37,\theta_{13}
\right)
-\textfrac{1}{8}\,\theta_3 \wedge \left(
13\,\theta_{12}
+96\,U\,\theta_{23}
+87\,U\,\xi^2
+45\,U^2\,\xi^3
\right)
\nonumber
\\ \fl
&-&
\textfrac{5}{8}\,\theta_{11}\wedge (\theta_{22}+\xi^2)
-2\,U\,\theta_{12} \wedge (\xi^2-U\xi^3)
+\theta_{13} \wedge \left(7\,\theta_{23}
-\textfrac{29}{8}\,U\,\xi^3\right),
\nonumber
\\ \fl
d\eta_7
&=&
\eta_8 \wedge \theta_0
+\eta_{11} \wedge \xi^1
+\eta_9 \wedge \xi^2
+\eta_{10} \wedge \xi^3
+\eta_3 \wedge (24\,U\,\theta_1-11\,\theta_{11})
-\eta_4 \wedge (3\,U\,\theta_0+2\,\theta_1)
\nonumber
\\ \fl
&-&
\textfrac{1}{4}\,\eta_2 \wedge \left(
13\,U^2\,\theta_0
-6\,U\,\theta_1
+4\,\theta_{11}
+8\,U\,\theta_{13}
\right)
-\textfrac{1}{8}\,\eta_5 \wedge \left(
19\,\theta_0 +16\,(\theta_3+U\,\xi^2)
\right)
\nonumber
\\ \fl
&-&
2\,\eta_6 \wedge (3\,\eta_2+2\,\theta_2+4\,U\,\xi^3)
-\textfrac{1}{8}\,\eta_7 \wedge \left(
20\,\eta_1 
-36\,(\theta_{22}+\xi^2)
-21\,\xi^3
\right)
%%%%
\nonumber
\\ \fl
&-&
\textfrac{1}{64}\,\theta_0 \wedge \left(
5\,(192\,U-1)\,\theta_1 
-88\,U^2\,\theta_2 
-32\,U\,\theta_3 
-440\,\theta_{11}
-272\,U\,\theta_{12}
-168\,\theta_{13} 
\right.
\nonumber
\\ \fl
&+&
\left.
16\,U^2\,\xi^2
-120\,U^3\,\xi^3
\right)
-\textfrac{1}{8}\,\theta_1 \wedge \left(
12\,\theta_3
+18\,\theta_{12}
+248\,U\,\theta_{23}
+224\,U\,\xi^2
-125\,U^2\,\xi^3
\right)
\nonumber
\\ \fl
&+&
\theta_2 \wedge \left(
36\,U\,\theta_1
-2\,U^2\,\theta_3
-\textfrac{109}{8}\,\theta_{11}
+3\,U\,\theta_{13} 
\right)
-\theta_3 \wedge (3\,U\,\theta_{12}+2\,\theta_{13}-U^2\,\xi^2)
\nonumber
\\ \fl
&+&
\theta_{11} \wedge \left(
11\,\theta_{23}
+\textfrac{77}{8}\,\xi^2 
-5\,U\,\xi^3
\right)
-3\,U\,\theta_{13} \wedge (\xi^2-\,U\,\xi^3)
\nonumber
\end{eqnarray}
}


\section*{Bibliography}
\begin{thebibliography}{99}
\bibitem{Blaszak} B{\l}aszak, M.: Classical R-matrices on Poisson algebras and related
    dispersionless systems, Phys. Lett. A {\bf 297}, 191--195 (2002)
\bibitem{BogdanovKonopelchenko} Bogdanov, L.V., Konopelchenko, B.G.: Nonlinear Beltrami
    equations and $\tau$-functions for dispersionless hierarchies, Phys. Lett. A {\bf 322},
    330--337 (2004) 
\bibitem{BCGGG} Bryant, R.L., Chern, S.S., Gardner, R.B., Goldschmidt, H.L., Griffiths, P.A.:
    Exterior Differential Systems. N.Y., Springer-Verlag (1991)
\bibitem{BryantGriffiths} Bryant, R.L., Griffiths, Ph.A.:
    Characteristic cohomology of differential systems (II):
    conservation laws for a class of parabolic equations,
    Duke Math. J. {\bf 78}, 531--676  (1995)
\bibitem{Cartan1}  Cartan, \'E.: Sur la structure des groupes infinis de transformations.
    {\OE}uvres Compl{\`e}tes,   Part II, {\bf 2}, p. 571--715.
    Gauthier - Villars, Paris (1953)
\bibitem{Cartan2}  Cartan, \'E.: Les sous-groupes des groupes continus de transformations.
    {\OE}uvres Compl{\`e}tes,   Part II, {\bf 2}, p. 719--856.
    Gauthier - Villars, Paris (1953)
\bibitem{Cartan3} Cartan, \'E.: Les probl\`emes d'\'equivalence. 
    {\OE}uvres Compl{\`e}tes,   Part II, {\bf 2}, p. 1311--1334.
    Gauthier - Villars, Paris (1953)
\bibitem{Cartan4} Cartan, \'E.: La structure des groupes infinis. 
    {\OE}uvres Compl{\`e}tes,   Part II, {\bf 2}, p. 1335--1384.
    Gauthier - Villars, Paris (1953)
\bibitem{ChangTu} Chang, J.-H., Tu, M.-H.: On the Miura map between the dispersionless KP and
    dispersionless modified KP hierarchies.  J. Math. Phys., {\bf 41}, 5391--5406 (2000).
\bibitem{DoddFordy} Dodd, R., Fordy, A.: The prolongation structures of quasipolynomial
    flows. Proc. Roy. Soc. London, A, {\bf 385}, 389--429  (1983)
\bibitem{Dunajski} Dunajski, M.: A class of Einstein--Weil spaces associated to an integrable
   system of hydrodynamic type, J. Geom. Phys. {\bf 51}, 126-137 (2004)  
\bibitem{Estabrook} Estabrook, F.B.: Moving frames and prolongation algebras.
    J. Math. Phys., {\bf 23} (1982), 2071--2076 (1982)
\bibitem{FelsOlver} Fels, M., Olver, P.J.: Moving coframes. I.
    A practical algorithm. Acta. Appl. Math. {\bf 51}, 161--213 (1998)
\bibitem{FerapontovKhusnutdinova} Ferapontov, E.V., Khusnutdinova, K.R.: The characterization of
    two-component (2+1)-di\-men\-si\-o\-nal integrable systems of hydrodynamic type, J. Phys.
     A.: Math. Gen. {\bf 37} 2949--2963 (2004) 
\bibitem{Gardner} Gardner, R.B.: The method of equivalence and its applications.
    CBMS--NSF regional conference series in applied math., SIAM, Philadelphia (1989)
\bibitem{Harrison1995}  Harrison, B.K.: On methods of finding B\"acklund transformations in
     systems with more than two independent variables. J. Nonlinear Math. Phys., {\bf 2}, 
     201--215 (1995)
\bibitem{Harrison2000}  Harrison, B.K.: Matrix methods of searching for Lax pairs and a paper by
    Est\'evez. Proc. Inst. Math. NAS Ukraine, {\bf 30}, Part 1, 17--24  (2000)
\bibitem{Hoenselaers} Hoenselaers, C.: More prolongation structures. 
    Prog. Theor. Phys. {\bf 75}, 1014--1029  (1986)
\bibitem{Igonin} Igonin, S.: Coverings and the fundamental group for partial differential
    equations. J. Geom. Phys., {\bf 56}, 939--998 (2006).
\bibitem{Kamran} Kamran, N.: Contributions to the Study of the Equivalence
    Pro\-blem of \'Elie Cartan and its Applications to Partial and
    Or\-di\-na\-ry Differential Equations.
    Mem. Cl. Sci. Acad. Roy. Belg., {\bf 45}, Fac. 7 (1989) 
\bibitem{KonopelchenkoAlonso} Konopelchenko, B., Mart\'inez Alonso, L.: Dispersionless scalar
    hierarchies, Whitham hierarchy and the quasi-classical $\bar\partial$-method, J. Math. Phys.
   {\bf43} 3807--3823 (2003)
\bibitem{KV84} Krasil'shchik, I.S., Vinogradov, A.M.: Nonlocal symmetries and the theory of
    coverings. Acta Appl. Math., {\bf 2}, 79--86  (1984)
\bibitem{KLV} Krasil'shchik, I.S., Lychagin, V.V., Vinogradov, A.M.:
    Geometry of Jet Spaces and Nonlinear Partial Differential Equations. Gordon and Breach, 
    New York (1986)  
\bibitem{KV89} Krasil'shchik, I.S., Vinogradov, A.M.: Nonlocal trends in the geometry of
    dif\-fe\-ren\-ti\-al equations: symmetries, con\-ser\-va\-ti\-on laws, and B\"acklund 
   trans\-for\-ma\-ti\-ons. Acta Appl. Math., {\bf 15}, 161--209  (1989)
\bibitem{KV99} Krasil'shchik, I.S., Vinogradov, A.M. (eds.):  Symmetries and
    Con\-ser\-va\-ti\-on Laws for Differential Equations of Ma\-the\-ma\-ti\-cal Physics.
    Transl. Math. Mo\-no\-graphs 182, Amer. Math. Soc., Pro\-vi\-dence (1999).
\bibitem{Krichever} Krichever, I.M.: The averaging method for two-dimensional "integrable"
    equations. Funct.Anal. Appl. {\bf 22}, 200--213 (1988)  
\bibitem{Kupershmidt} Kupershmidt, B.A.: The quasiclassical limit of the modified KP hierarchy.
    J. Phys. A Math. Gen. {\bf 23}, 871--886 (1990)
\bibitem{Kuzmina} Kuz'mina, G.M.:  On a possibility to reduce a system of two first-order
    partial differential equations to a single equation of the second order. Proc. Moscow State
    Pe\-da\-gog. Inst. {\bf 271}, 67--76  (1967) (in Russian)
\bibitem{Marvan1992} Marvan, M.: On zero-curvature representations of partial differential
     equations.  Proc. Conf. on Diff. Geom. and Its Appl., Opava (Czech Republic),  103--122
     (1992)
\bibitem{Marvan1997} Marvan, M.:  A direct procedure to compute zero-curvature representations.
    The case $\mathfrak{sl}_2$. In: Proc. Int. Conf. on Secondary
    Calculus and Cohomological Physics, Mos\-cow, Russia, August 24-31, 1997. Available via the
    Internet at ELibEMS, 
    {\tt http://www.emis.de/proceedings.}
\bibitem{Marvan2002} Marvan, M.: On the horizontal gauge cohomology and nonremovability of 
     the spectral parameter. Acta Appl. Math {\bf 72}, 51--65 (2002)
\bibitem{Morozov2002} Morozov, O.I.: Moving coframes and symmetries of
    differential equations. J. Phys. A, Math. Gen., {\bf 35}, 2965--2977  (2002)
\bibitem{Morozov2006} Morozov, O.I.: Contact-equivalence problem for linear
    hyperbolic equations.  J. Math. Sci., {\bf 135}, 2680--2694  (2006)
\bibitem{Morozov2007}  Morozov, O.I.: Coverings of differential equations and Cartan's
    structure theory of Lie pseudo-groups.  Acta Appl. Math. {\bf 99}, 309--319  (2007)
\bibitem{Morozov2008} Morozov, O.I.:  Cartan's structure theory of symmetry pseudo-groups,
    coverings and multi-valued solutions for the  Khokhlov--Zabolotskaya equation, 
    Acta Appl. Math. {\bf 101}, 231-241   (2008)
\bibitem{Morozov2008b} Morozov, O.I.: Cartan's structure of symmetry pseudo-group and coverings 
    for the r-th modified dispersionless Kadomtsev--Petviashvili equation. Preprint 
    {\tt arXiv:0803.096v1 [nlin.SI]}  (2008)
\bibitem{Morris1976} Morris, H.C.: Prolongation structures and nonlinear evolution equations 
    in two spatial  dimensions. J. Math. Phys., {\bf 17}, 1870--1872  (1976)
\bibitem{Morris1979} Morris, H.C.: Prolongation structures and nonlinear evolution equations 
    in two spatial  dimensions: a general class of equations. J. Phys. A, Math. Gen., {\bf 12}, 
    261--267 (1979)
\bibitem{Olver95} Olver, P.J.: Equivalence, Invariants, and Symmetry. Cambridge, Cambridge
    Uni\-ver\-si\-ty Press (1995)
\bibitem{Sakovich} Sakovich, S.Yu.: On zero-curvature representations of evolution equations. 
    J. Phys. A, Math. Gen., {\bf 28}, 2861--2869  (1995)
\bibitem{Stormark} Stormark, O.: Lie's Structural Approach to PDE Systems. Cambridge, Cambridge
    Uni\-ver\-si\-ty Press (2000)
\bibitem{Takasaki} Takasaki, K.: Quasi-classical limit of BKP hierarchy and W-infinity 
    symmetries, Lett. Math. Phys. {\bf 28}, 177--185 (1993)
\bibitem{Tondo} Tondo, G.S.: The eigenvalue problem for the three-wave resonant interaction in
   (2+1)  dimensions via the prolongation structure. Lett. Nuovo Cimento, {\bf 44}, 297--302 
   (1985)
\bibitem{Vasil'eva} Vasil'eva, M.V.: The Structure of Infinite Lie Groups of Transformations.
    Moscow, MGPI (1972) (in Russian)
\bibitem{WE} Wahlquist, H.D., Estabrook F.B.: Prolongation structures of nonlinear 
    evo\-lu\-ti\-on  equations.  J. Math. Phys., {\bf 16}, 1--7  (1975)
\bibitem{Zakharov82} Zakharov, V.E.: Integrable systems in multidimensional spaces. 
    Lect. Notes Phys., {\bf 153}, 190--216  (1982)
\end{thebibliography}
\end{document}